\input amstex
\documentstyle{amsppt}
\input epsf 
\pagewidth{6truein}
\pageheight{8.8truein}
\NoBlackBoxes
\def\A{{\Cal A}}
\def\D{{\Cal D}}
\def\F{{\Cal F}}
\def\G{{\Cal G}}
\def\K{{\Cal K}}
\def\L{{\Cal L}}
\def\M{{\Cal M}}
\def\N{{\Cal N}}
\def\Z{{\Cal Z}}
\def\SS{{\Cal S}}
\def\Sp{{\Cal S}p}
\def\cee{{\Bbb C}}
\def\nat{{\Bbb N}}
\def\Ca{{\text{\bf Ca}}}
\def\bK{{\text{\bf K}}}
\def\bT{{\text{\bf T}}}
\def\bX{{\text{\bf X}}}
\def\bY{{\text{\bf Y}}}
\def\bZ{{\text{\bf Z}}}
\def\bn{{\text{\bf n}}}
\def\cb{\operatorname{cb}}
\def\Ex{\operatorname{Ex}}
\def\Id{\operatorname{Id}}
\def\MAX{\operatorname{MAX}}
\def\MIN{\operatorname{MIN}}
\def\NEW{\operatorname{NEW}}
\def\op{\operatorname{op}}
\def\ep{\varepsilon}
\def\bPi{{\boldsymbol \Pi}}
\def\To{\Rightarrow}
\def\chix{{\raise.5ex\hbox{$\chi$}}}
\def\defeq{\mathop{\buildrel {\text{df}}\over =}}
\def\faceq{\mathop{\buildrel {\text{Fact}}\over =}}
\def\Ba{\mathop{Ba}\nolimits}
\def\Sp{\mathop{{\Cal S}p}\nolimits}
\def\oQ{\overline{Q}}
 
\topmatter
\title On certain extension properties for the space of compact operators
\endtitle
\rightheadtext{THE SPACE OF COMPACT OPERATORS}
\author Timur Oikhberg and Haskell P. Rosenthal\endauthor
\address Timur Oikhberg, Department of Mathematics, The University 
of Texas at Austin, Austin, TX 78712-1082\endaddress
\email timur\@math.utexas.edu\endemail
\address Haskell Rosenthal, Department of Mathematics, The University
of Texas at Austin, Austin, TX 78712-1082\endaddress
\email rosenthl\@math.utexas.edu\endemail
\abstract
Let $Z$ be a fixed separable operator space, $X\subset Y$ general 
separable operator spaces, and $T:X\to Z$ a completely bounded map. 
$Z$ is said to have the Complete Separable Extension Property (CSEP) if 
every such map admits a completely bounded extension to $Y$; the 
Mixed Separable Extension Property (MSEP) if every such $T$ admits a bounded 
extension to $Y$. 
Finally, $Z$ is said to have the Complete Separable Complementation Property 
(CSCP) if $Z$ is locally reflexive and $T$ admits a completely bounded 
extension to $Y$ {\it provided\/} $Y$ is locally reflexive and $T$ is a 
complete surjective isomorphism. 
Let $\bK$ denote the space of compact operators on separable Hilbert space 
and $\bK_0$ the $c_0$ sum of $\M_n$'s (the space of 
``small compact operators''). 
It is proved that $\bK$ has the CSCP, using the second author's previous 
result that $\bK_0$ has this property. 
A new proof is given for the result (due to E.~Kirchberg) that $\bK_0$ 
(and hence $\bK$) fails the CSEP. 
It remains an open question if $\bK$ has the MSEP; it is proved this is 
equivalent to whether $\bK_0$ has this property. 
A new Banach space concept, Extendable Local Reflexivity (ELR), is 
introduced to study this problem. 
Further complements and open problems are discussed.
\endabstract
\endtopmatter

\document

\head Contents\endhead 

{\narrower\smallskip

\line{\qquad Introduction\dotfill 2}
\line{\S1. Extending Complete Isomorphisms into $B(H)$\dotfill 4}
\line{\S2. An operator space construction on certain subspaces of $\M_\infty$
\dotfill 11}
\line{\S3. The $\lambda$-Mixed Separable Extension Property and
Extendably Locally\hfill}
\line{\qquad Reflexive Banach spaces\dotfill 17}
\line{\S4. $\bK_0$ fails the CSEP: a new proof and generalizations
\dotfill 33}
\line{\qquad References\dotfill 40}

}

\newpage
\head Introduction\endhead 

The space $\bK$ of compact operators on a separable infinite dimension 
Hilbert space $H$ is often that thought of as the non-commutative analogue 
of $c_0$, the space of sequences vanishing at infinity. 
Indeed, if one regards $\bK$ as matrices with respect to a fixed 
orthonormal basis of $H$, the {\it diagonal\/} matrices form a subalgebra 
isometric to $c_0$. 
In 1941, A.~Sobcyk proved that $c_0$  has the Separable Extension Property 
(SEP) \cite{S}: 
If $Z=c_0$, then given $X\subset Y$ separable Banach spaces and 
$T:X\to Z$ a bounded linear operator, there exists a bounded linear 
operator $\tilde T:\to Z$ extending $T$. 
In 1977, M.~Zippin proved the (much deeper!) converse to this result 
\cite{Z}, 
{\it any infinite-dimensional separable Banach space $Z$ with the {\rm SEP} 
is isomorphic to $c_0$.}
We continue here the study of operator space analogues of the SEP, 
initiated in \cite{Ro2}, with the goal in particular of specifying which of 
these analogues $\bK$ satisfies. 
(For basic facts about operator spaces see \cite{Pi3}; also see the 
Introduction to \cite{Ro2} for a brief summary and orientation.) 

Thus we consider a fixed operator space $Z$, and consider the following 
diagram:
$$\lower.5truein\hbox{\epsfysize=0.8truein\epsfbox{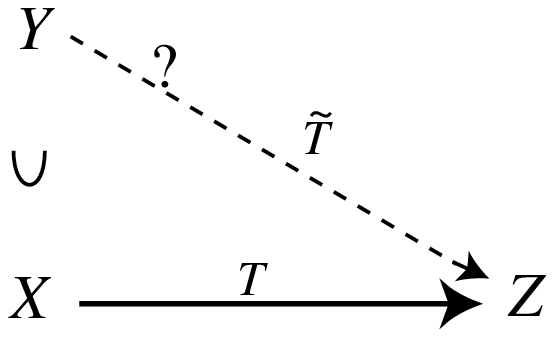}}$$
Here, $X$ and $Y$ are (appropriately general) separable operator spaces 
and $T$ is a completely bounded linear map.

$Z$ is said to have the Complete Separable Extension Property (CSEP) 
if every such $T$ admits a 
{\it completely bounded linear extension\/} $T$; 
the Mixed Separable Extension Property (MSEP) if $T$ admits a 
{\it bounded linear extension\/} $\tilde T$, and the 
Complete Separable Complementation Property (CSCP) if $T$ admits a 
bounded linear extension $\tilde T$ {\it provided\/} $Z$ is separable 
locally reflexive, $Y$ is also locally reflexive, {\it and\/} $T$ is a 
complete surjective isomorphism. 
If $1\le \lambda$ is such that $\tilde T$ can be chosen with 
$\|\tilde T\|_{\cb} \le\lambda \|T\|_{\cb}$ in the CSEP-case, we say $Z$ 
has the $\lambda$-CSEP; if $\|\tilde T\|\le\lambda\|T\|_{\cb}$ in 
the MSEP-case, we say $Z$ has the $\lambda$-MSEP. 
It follows easily that if $Z$ {\it has\/} the CSEP (resp. the MSEP), then 
$X$ has the $\lambda$-CSEP (resp. the $\lambda$-MSEP) for some $\lambda\ge1$.

Of course these properties are intimately connected with injectivity 
notions; thus $Z$ is called (isomorphically) injective (resp. mixed 
injective) if this diagram admits a completely bounded solution (resp. 
bounded solution) $T$ for arbitrary (not necessarily separable) 
operator spaces $X$ and $Y$. 
As in the separable setting, if $Z$ is injective (resp. mixed injective), 
there is a $\lambda\ge1$ so that $\tilde T$ may always be chosen with 
$\|\tilde T\|_{\cb} \le \lambda\|T\|_{\cb}$ (resp. $\|\tilde T\|\le 
\lambda \|T\|_{\cb}$); if $\lambda$ works, we say $Z$ is $\lambda$-injective 
(resp. $\lambda$-mixed injective). 
We say $Z$ is isometrically injective (resp. isometrically mixed injective) 
when $\lambda=1$. 

It is a fundamental theorem in operator space theory that $B(H)$ is 
isometrically injective for any Hilbert space $H$, where $B(H)$ denotes 
the space of bounded linear operators on $H$. 
It follows easily that if $X$ is an operator space with $X\subset B(H)$ 
for some Hilbert space $H$, then  $X$ is isomorphically injective 
(resp. mixed injective) if and only if $X$ is completely complemented 
(resp. complemented) in $B(H)$. 

The separable extension properties we consider have their primary interest 
for $\lambda\ge2$. 
Indeed, if $Z$ is separable, 
then if $\lambda<2$ and $Z$ has the $\lambda$-CSEP, 
it is proved in \cite{Ro2} that $Z$ is $\lambda$-injective; we show 
analogously here that if $Z$ has the $\lambda$-MSEP, $Z$ is $\lambda$-mixed 
injective (and moreover $Z$ is reflexive, whence by a result of G.~Pisier, 
$Z$ is actually Hilbertian (cf. \cite{R})). 

One of the main results of this work is that $\bK$ has the CSCP. 
A result of E.~Kirchberg yields that $\bK$ {\it fails\/} the CSEP \cite{Ki1}. 
We give a new proof and further complements in Section~4. 

It is proved in \cite{Ro2} that $\bK_0$ has the CSCP, where $\bK_0$ 
denotes the space of ``small compact operators'', namely the $c_0$-sum 
of $\M_n$'s, where $\M_n$ denotes the space of complex $n\times n$ 
matrices, identified with $B(\cee^n)$ for all $n$, $\cee^n$ being the 
standard $n$-dimensional complex Hilbert space. 
We obtain that $\bK$ has the CSCP via the following route: 
in Section~1, we show that if $X\subset Y$ are given separable operator 
spaces, then {\it any complete isomorphism from $X$ into $B(H)$ admits 
a complete isomorphic extension from $Y$ into $B(H)$} (Theorem~1.1). 
It follows from this result that if $X\subset B(H)$ is fixed with $X$ 
separable locally reflexive, then {\it $X$ has the {\rm CSCP} provided $X$ 
is completely complemented in $Y$ for any separable locally reflexive 
operator space $Y$ when $X\subset Y\subset B(H)$} (see Corollary~1.8). 
Now it follows from the main result of section~2 (Theorem~2.1) that if 
$\bK \subset Y \subset B(H)$ (where this is the {\it natural\/} embedding 
of $\bK$ in $B(H)$) with $Y$ separable, there is an absolute constant $C$ 
and for all $\ep>0$, a projection $P$ on $B(H)$ with $\|P\|_{\cb}<1+\ep$ 
with $Y$ and $\bK$ invariant under $P$ so that $(I-P)Y\subset\bK$ and 
$d_{\cb} (\bK,\bK_0)\le C$. 
It then easily follows that $\bK$ is completely complemented in $Y$ 
provided $Y$ is locally reflexive, from the fact that then $P\bK$ has this 
property by the result in \cite{Ro2}. 
We do not know if $\bK$ has the MSEP. 
However Theorem~2.1 also yields that $\bK$ has the MSEP if $\bK_0$ has 
this property (Proposition~2.3). 

We also obtain in Section~1 that if an operator space $Z$ has the CSCP, it has 
the following stronger property: 
{\it there is a completely bounded operator $\tilde T$ completing the above 
diagram whenever $Y$ is separable 
locally reflexive and $X$ is locally complemented in $Y$} (Theorem~1.4). 
As shown in \cite{Ro2}, $X$ ``automatically'' is locally complemented 
provided $X$ is completely isomorphic to a nuclear $C^*$-algebra, or more 
generally, if $X^{**}$ is isomorphically injective. 
($X$ is called locally complemented in $Y$ provided there is a $C<\infty$ 
so that $X$ is $C$-completely complemented in $W$ for all $X\subset W\subset Y$
with $W/X$ finite-dimensional). 
It was also previously proved in \cite{Ro2} that $\bK_0$ has this 
stronger property, 
and moreover one may drop the assumption that $Y$ is locally reflexive. 

The MSEP is studied in Section~3, where we introduce the following 
concept: 
Given operator spaces $X$ and $Y$, $X$ is called {\it completely 
semi-isomorphic to\/} $Y$ if there is a completely bounded 
surjective map 
$T:X\to Y$ which is a Banach isomorphism; $X$ is called completely 
semi-isometric to $Y$ in case $T$ can be chosen with $\|T\|_{\cb} = 1=
\|T^{-1}\|$. 
We then have the simple permanence property: 
mixed injectivity and the MSEP are both preserved under complete 
semi-isomorphisms (Proposition~3.9). 
The finite-dimensional isometrically mixed injectives are known up to 
Banach isometry; they are the $\ell^\infty$-direct sums of Cartan factors 
of type~IV (see Theorem~A, following Problem~3.2). 
This result suggests a possible classification of the isometrically 
injective finite-dimensional operator spaces; are all such completely 
semi-isometric to an $\ell^\infty$-direct sum of Cartan factors of 
types~I--IV?  (Problem 3.3). 
A remarkable factorization result of M.~Junge's and the semi-isomorphism 
concept yield that the classification problem of the finite-dimensional 
mixed injectives is exactly analogous to the famous open commutative 
case; namely {\it if $X$ is finite-dimensional and $\lambda$-mixed 
injective, then for all $\ep>0$, there is 
an $n$ so that $X$  is $(\lambda+\ep)$-semi-isomorphic 
to some $(\lambda+\ep)$-complemented subspace of $\M_n$} (Proposition~3.10). 

To further penetrate the MSEP-problem for $\bK$, we introduce a new pure 
Banach space concept in Section~3, that of Extendable Local Reflexivity (ELR). 
Several equivalences are given in Proposition~3.12; Theorem~3.13 yields 
the remarkable equivalence that a Banach space $X$ is Extendably Locally 
Reflexive and has the bounded approximation property if and only if $X^*$ 
has the bounded approximation property. 
The ``if'' part of this assertion is quite simple, and was discovered by 
the second author of this present paper in the fall of 1997. 
The remarkable ``only if'' part is due to W.B.~Johnson and the first 
author of the present paper \cite{JO}. 
Actually, ELR has a complete analogue, and the complete version of 
Theorem~3.13 also holds. 
The motivation for the introduction of this concept: 
{\it if $X$ is a separable operator space so that $X^{**}$ is isomorphically 
mixed injective and is (Banach) {\rm ELR}, then $X$ has the {\rm MSEP}}  
(Theorem~3.14). 
In particular, $\bK$ has the MSEP if $B(H)$ is ELR. 
The proof of this result involves a construction mixing Banach and 
operator space ideas, perhaps of interest in its own right (Lemma~3.16). 

Section 4 establishes some necessary conditions for certain operator spaces 
to have the CSEP, yielding in particular a ``qualitative'' proof that 
$\bK_0$ (and hence $\bK$) fails the CSEP. 
It is proved for example that 
{\it if $Z_1,Z_2,\ldots$ are finite-dimensional operator spaces, then if 
$(Z_1\oplus Z_2\oplus\cdots)_{c_0}$ has the {\rm CSEP}, 
$\{Z_1,Z_2,\ldots\}$ is of finite matrix type, and the $Z_j$'s are all 
$\lambda$-injective for some $\lambda$} (Corollary~4.2). 
The converse to this result is established in Proposition~2.15 of 
\cite{Ro2}. 
(See the beginning of Section~4 for the definition of {\it finite matrix
type\/}.) 
It is further shown that if $Z$ is a separable operator space so that 
$c_0(Z)$ has the CSEP, then $Z$ is of finite matrix type (Corollary~4.4). 
We conjecture that if $Z$ is separable with the CSEP, then $Z$ itself 
is of finite matrix type. 

Corollary 4.2 yields that $\bK_0$ fails the CSEP, for $\{\M_1,\M_2,\ldots\}$ 
is {\it not\/} of finite matrix type. 
We obtain a stronger quantitative result in Corollary~4.9: 
{\it For every $n$, there exists an operator space $Y_n$ containing 
$\bK_0$ so that $Y_n/\bK_0$ is completely isometric to $\ell_n^\infty$, 
yet $\bK_0$ is not $\lambda$-completely complemented in $Y_n$ if\/} 
$\lambda\le \sqrt{n}/2$. 
It follows then from results in \cite{Ro2} that $Y_n$, which is of course 
locally reflexive, is not $\lambda$-locally reflexive if $\lambda < 
(\sqrt{n}/2)-3$. 
Putting these $Y_n$'s together, we then obtain an operator space $Y$ 
containing $\bK_0$ so that $\bK_0$ 
{\it is not completely complemented in $Y$, yet $Y/\bK_0$ is completely 
isometric to\/} $c_0$. 
Thus $Y$ cannot be locally reflexive since $\bK_0$ has the CSCP 
(of course this also follows by its construction). 
It also then follows by results of E.~Kirchberg (\cite{Ki2}) that $Y$ 
is not a nuclear operator space; however $\bK_0$ and $Y/\bK_0$ are 
obviously nuclear. 
We also show in Proposition~4.11 that any descending sequence of 1-exact
finite dimensional Banach isometric spaces must be bounded below. 

We finally show that $Z= \bK_0$ fails to have a completely bounded 
solution $\tilde T$ to the above diagram if $Y$ is general locally 
reflexive separable, $X$ a general subspace. 
Actually, we obtain that there exists a subspace $X$ of $C_1$ (the operator 
space of trace class operators) and a completely bounded linear map 
$T:X\to \bK_0$ so that $T$ has no completely bounded extension $\tilde T$ 
to $C_1$. 
(A remarkable result due to M.~Junge yields that $C_1$ is 1-locally reflexive.) 
In fact, we establish in Proposition~4.14 that 
{\it if $Z$ is separable and there is a completely bounded solution to 
the above diagram for arbitrary $X\subset C_1$, $Y=C_1$, then $Z$ has the\/} 
CSEP.

\head {\smc Section} 1\\
{\bf Extending complete isomorphisms into $B(H)$}\endhead 

The main result of this section is the following:

\proclaim{Theorem 1.1} 
Let $Y$ be a separable operator space, $X$ a subspace of $Y$, and 
$T:X\to B(H)$ a complete isomorphic injection of $X$. 
There exists a complete isomorphic injection $\tilde T:Y\to B(H)$ 
extending $T$.
\endproclaim

\remark{Remarks}
1. We obtain 
$$\|\tilde T\|_{\cb} \le 3\|T\|_{\cb} \text{ and } 
\|\tilde T\|_{\cb} \|\tilde T^{-1}\|_{\cb} 
\le 12\|T\|_{\cb} \|T^{-1}\|_{\cb} +6\ .
\tag $*$ $$

2. Our proof of the Theorem uses ideas from \cite{LR}. 
In face, our argument may be refined to obtain the following stronger 
result, analogous to a result in \cite{LR}, showing that completely 
isomorphic separable sequences of $B(H)$ lie in the same position. 
{\it Let $X,X'$ be separable operator subspaces of $B(H)$ and $T:X\to X'$ 
a complete surjective isomorphism. 
There exists a complete surjective isomorphism $\tilde T:B(H)\to B(H)$ 
extending $T$.} 

We first give an operator-space version of a result of A.~Pe{\l}czy\'nski 
\cite{Pe1}, for which we use the following lemma (which is quite different 
than the argument in \cite{Pe1}). 
\endremark

\proclaim{Lemma 1.2}
Let $X\subset Y$ and $Z$ be Banach spaces and $T:Y\to Z$ a bounded linear 
operator so that $T|X$ is an (into) isomorphism. 
Let $\Pi :Y\to Y/X$ be the quotient map and define $\tilde T:Y\to 
Z\oplus Y/X$ by 
$$\tilde Ty = Ty \oplus \Pi y\ \text{ for all } y\in Y\ .
\tag 1$$
Then $\tilde T$ is an into-isomorphism with $\tilde T|X=T$. 
In fact 
$$\|\tilde T\|\, \|\tilde T^{-1}\| \le 2\|T\|\, \|(T|X)^{-1}\| + 1
\tag 2$$
\endproclaim

\remark{Remark}
We put the $\infty$-norm on the direct sum $Z\oplus Y/X$. 
If $\tilde X= Tx$, $\tilde Y = Ty$, $(T|X)^{-1}$ refers to the inverse 
map from $\tilde X$ to $X$, $T^{-1}$ the map corresponding inverse map 
from $\tilde Y$ to $Y$.
\endremark

\demo{Proof}
It is trivial that $\tilde T|X=Z$. 
(Of course we identity $Z$ with $Z\oplus 0$.) 
We may assume without loss of generality that $\|T\|=1$. 
Let $\delta  = \|(T|X)^{-1}\|^{-1}$, and fix $y\in Y$ with $\|y\|=1$. 
Set $\tau = \|\Pi y\|$; let $\ep >0$, and choose $x\in X$ with 
$\|x-y\| \le \tau +\ep$. 
Then 
$$\|x\| \ge 1-\tau\ .
\tag 3$$
Hence 
$$\|Tx\| \ge \delta (1-\tau)\ ,
\tag 4$$
and so 
$$\|\tilde Ty\| \ge \delta (1-\tau)- (\tau +\ep)\ .
\tag 5$$
Of course also 
$$\|\tilde T y\| \ge \|\Pi y\| = \tau\ .
\tag 6$$
Hence 
$$\align 
\|\tilde Ty \| & \ge \max \{\delta - \ep - (\delta +1)\tau,\tau\}\\
&\ge \frac{\delta-\ep}{\delta+2}\ \text{ for }\ 0\le \tau\le1\ .
\endalign$$
Since $\ep>0$, we have proved that 
$$\| \tilde T^{-1}\| \le 2\|(T|X)^{-1}\| +1
\tag 7$$
which establishes (2) and thus the Lemma.\qed
\enddemo

The next result yields \cite{Pe1, Proposition 1} when restricted to 
the Banach space category. 

\proclaim{Proposition 1.3} 
Let $X\subset Y$, $\tilde X$ be operator spaces and $T:X\to \tilde X$ 
a complete surjective isomorphism. 
There exists an operator space $\tilde Y\supset \tilde X$ and a complete 
surjective isomorphism $\tilde T:Y\to \tilde Y$ extending $T$, in fact    
satisfying 
$$\|\tilde T\|_{\cb} = \|T\|_{\cb}\text{ and } 
\|\tilde T\|_{\cb} \|\tilde T^{-1}\|_{\cb} \le 2 \|T\|_{\cb} 
\|T^{-1}\|_{\cb} +1\ .
\tag 8$$
\endproclaim

\remark{Remark}
Proposition 1.3 (or rather its proof) is used in the proof of the main 
result of this section, Theorem~1.1. 
However if we use an alternate construction, due to G.~Pisier (see 
10b, p.137 of \cite{Pi2}), we obtain considerably better quantitative 
information. 
We first formulate the result, then give the construction, leaving the 
details of proof to the interested reader.
\endremark

\proclaim{Proposition}
Let $X\subset Y$, $\tilde X$ be operator spaces and $T:X\to \tilde X$ a 
completely bounded map. 
There exists an operator space $\tilde Y\supset \tilde X$ and a completely 
bounded map $\tilde T:Y\to \tilde Y$ extending $T$, satisfying 
the following: 
\roster
\item"(i)" $\|\tilde T\|_{\cb} = \|T\|_{\cb}$.
\item"(ii)" If $T$ is an (into) isomorphism (resp. complete isomorphism) so 
is $\tilde T$, and $\|\tilde T^{-1}\| = \|T^{-1}\|$ (resp. $\|\tilde T^{-1}
\|_{\cb} = \|T^{-1}\|_{\cb}$).
\item"(iii)" If $T$ is surjective, so is $\tilde T$.
\item"(iv)" $Y/X$ is completely isometric to $\tilde Y/\tilde X$.
\endroster
\endproclaim 

\proclaim{Corollary} 
Let $X,Y,\tilde X$ be operator spaces with $X\subset Y$ and $\tilde X$ 
and $X$ completely isomorphic; set $\beta = d_{\cb}(\tilde X,X)$. 
Then given $\ep>0$, there exists an operator space $\tilde Y\supset 
\tilde X$ and a complete bijection $T:Y\to \tilde Y$ with $TX=\tilde X$
and $\|T\|_{\cb} \|T^{-1}\|_{\cb} < \beta +\ep$. 
\endproclaim 

Pisier's construction (adapted to this setting) goes as follows: 
Let $Z = (Y\oplus \tilde X)_1$ (endowed with its natural operator space 
structure, cf. \cite{Pi3}). 
Assume first that $\|T\|_{\cb}=1$. 
Let $\Gamma = \{x\oplus -Tx:x\in X\}$. 
Then $\Gamma$ is a closed subspace of $Z$. 
Let $\tilde Y = Z/\Gamma$ and $\Pi :Z\to \tilde Y$ be the quotient map.
Define $j:\tilde X\to \tilde Y$ and $\tilde T:Y\to \tilde Y$ by 
$j(\tilde x) = \Pi (0\oplus \tilde x)$ and $\tilde T(y) = \Pi (y\oplus 0)$ 
for all $\tilde x\in\tilde X$ and $y\in Y$.

One then has the commutative diagram
$$\CD
Y @>{\tilde T}>> \tilde Y\\
@AiAA @AAjA\\
X @>T>> \tilde X
\endCD$$
where $i:X\to Y$ is the inclusion map. 
Then it follows that $j$ is a {\it complete $($into$)$ isometry\/} 
and one now verifies all the details of the Proposition (identifying 
$\tilde X$ with $j(\tilde x)$). 
In general, assuming $T\ne 0$, simply set $U= T/\|T\|_{\cb}$  and now 
carry out the construction for $U$ instead, but then simply define 
$\tilde T$ by $\tilde T(y) = \|T\|_{\cb} \Pi (y\otimes 0)$ for all $y\in Y$. 
E.g., to verify the first part of (ii), let 
$\delta = 1/\|T^{-1}\|$ (resp. $1/\|T^{-1}\|_{\cb}$) (and assume still 
$\|T\|_{\cb} =1$, so $\delta \le 1$); 
if $y\in Y$,
$$\eqalign{
\|\tilde T y\| &= \inf \{ \|y-x\| + \|Tx\| :Y\in X\}\cr
&\ge \inf \{ \delta \|y-x\| + \delta \|x\| :x\in X\}\cr
&\ge \delta \|y\|\ .\cr}$$
Hence also $\|\tilde T\|^{-1} = \frac1{\delta} = \|T^{-1}\|$. 

\demo{Proof of Proposition 1.3} 
We may assume $\tilde X\subset B(H)$ for a suitable Hilbert space $H$; 
also we may assume that $\|T\|_{\cb} =1$. 
Using the isometric injectivity of $B(H)$, choose $T':Y\to B(H)$ a linear 
map extending $T$ with also $\|T'\|_{\cb}=1$. 
Now we apply the result in Lemma~1.2 to $Z= B(H)$. 
Let $\Pi :Y\to Y/X$ be the quotient map and define 
$\tilde T:Y\to B(H) \oplus Y/X$ by 
$$\tilde Ty = T'y \oplus \Pi y\ \text{ for all }\ y\in Y\ .
\tag 9$$
Now setting $\tilde Y = \tilde T(Y)$, it follows from Lemma~1.2 that 
$\tilde Y$ is a closed linear subspace of $B(H)\oplus Y/X$ and $\tilde T$ 
is an isomorphism from $Y$ onto $\tilde Y$. 
We claim further that $\tilde T$ is indeed a complete isomorphism. 
Now it follows immediately from (9) and the complete contactivity of 
$T'$ that $\|\tilde T\|_{\cb}=1$. 
Let $\bX = \bK \otimes_{\op} X$, 
$\bY = \bK \otimes_{\op} Y$, and 
$\bT = I\otimes \tilde T :\bY \to \bK \otimes_{\op} B(H) \defeq \bZ$
where $I = \Id |\bK$. 
It follows that $\bT|\bX$ is an isomorphism onto $\bK \otimes_{\op}
\tilde X$ and 
$$\|(\bT|\bX)^{-1}\| = \|T^{-1}\|_{\cb}, \|\bT\| = \|\tilde T\|_{\cb}=1\ .
\tag 10$$
Thus we may apply Lemma 1.2. 
Let  $\bPi :\bY \to \bY/\bX$ be the quotient map and define  $\tilde{\bT}$ 
as in (1). 
Thus  Lemma~1.2 yields that $\tilde{\bT}$ is an isomorphism onto 
$\bK\otimes_{\op} \tilde Y$ and (2) yields
$$\|(\tilde{\bT})^{-1}\| \le 2\|(\tilde{\bT}|\bX)^{-1}\| +1\ .
\tag 11$$ 
But also 
$$\|(\tilde{\bT})^{-1}\| = \|(\tilde T)^{-1}\|_{\cb}\ .
\tag 12$$ 
(10)--(12) thus complete the proof.\qed
\enddemo

The next result is a simple application of Proposition~1.3, 
following the concept introduced in \cite{Ro2}: 
{\it A locally reflexive separable operator space $Z$ is said to have the 
Complete Separable Complementation Property {\rm (CSCP)} provided every 
complete isomorph of $Z$ is completely complemented in every separable 
locally reflexive operator superspace.}
Equivalently, given separable operator spaces $X\subset Y$ with $Y$ locally 
reflexive and $T:X\to Z$ a complete surjective isomorphism, there exists 
a completely bounded $\tilde T:Y\to Z$ extending $T$. 
That is, we have the diagram
$$\lower.5truein\hbox{\epsfysize=0.8truein\epsfbox{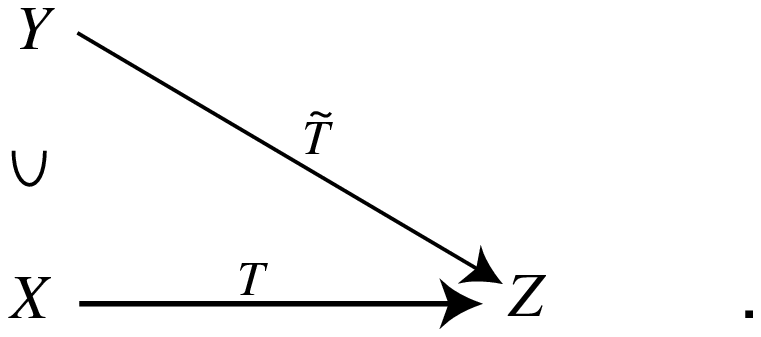}}
\tag 13$$

Our next result yields some equivalences for the CSCP 
We obtain a considerably stronger version of one of the 
equivalences at the end of this section (using Theorem~1.1). 

\proclaim{Theorem 1.4} 
Let $X$ be a separable locally reflexive operator space. 
Then the following are equivalent:
\roster
\item"(a)" $X$ has the CSCP.
\item"(b)" $X$ is completely complemented in every separable locally 
reflexive superspace. 
\item"(c)" $X$ is locally complemented in every separable locally reflexive 
operator superspace, and whenever $\tilde X$ is a locally complemented 
subspace of a locally reflexive separable operator superspace $\tilde Y$ and 
$T:\tilde X\to X$ is a completely bounded map, there exists a 
completely bounded map $\tilde T:\tilde Y\to X$ extending $T$. 
\endroster
\endproclaim

\remark{Remarks} 

1. If $X$ and $Y$ are operator spaces with $X\subset Y$, $X$ is said to be 
{\it locally complemented\/} in $Y$ if there is a $C$ so that $X$ is 
$C$-completely complemented in $Z$ for all linear subspaces $Z$ with 
$X\subset Z\subset Y$ and $Z/X$ finite dimensional. 
(Then one says $X$ is $C$-locally complemented in $Y$). 
A simple compactness argument yields that if $X$ is $(C+\ep)$-locally 
complemented in $Y$ for all $\ep>0$, $X^{**} = X^{\bot\bot}$ is 
$C$-completely complemented in $Y^{**}$ via a weak*-continuous projection. 
Conversely, it is proved in Sublemma~3.11 of \cite{Ro2} that 
{\it if $Y$ is locally reflexive and $X^{**}$ is completely complemented 
in $Y^{**}$, $X$ is locally complemented  in $Y$}. 

2. There  are many situations in which the hypotheses  for $\tilde X$ 
hold without $\tilde X$ necessarily being complemented. 
For example, if $\tilde X$ is completely isomorphic to a commutative 
$C^*$-algebra or more generally,  a nuclear $C^*$-algebra, then if 
$\tilde X\subset\tilde Y$, with $\tilde Y$ locally reflexive, $\tilde X$ 
is ``automatically'' locally complemented in $\tilde Y$, but e.g., in the 
commutative case, there are examples where $\tilde X$ and $\tilde Y$ are 
actually commutative separable $C^*$-algebras with $\tilde X$ uncomplemented 
in $\tilde Y$. 
Perhaps the most general hypothesis on $\tilde X$ alone:
{\it if $(\tilde X)^{**}$ is an isomorphically injective operator space, then 
$\tilde X$ is locally complemented in any locally reflexive operator 
superspace\/} (cf. \cite{Ro2}).

\demo{Proof} 
(a) $\To$ (b) is trivial. 

(b) $\To$ (a): Let $\tilde X$ and $\tilde Y$ be operator spaces with 
$\tilde X\subset \tilde Y$, $\tilde Y$ locally reflexive separable, 
and $\tilde X$ completely isomorphic to $X$. 
Let $T:\tilde X\to X$ be a complete surjective isomorphism. 
By Proposition~1.3, there exists an operator space $Y\supset X$ and a 
complete surjective isomorphism $\tilde T:\tilde Y\to Y$ extending $T$.
Then $Y$ is also separable and locally reflexive, hence there exists a 
completely bounded projection $P:Y\to X$. 
Then $\tilde P \defeq (T|X)^{-1} P\tilde T$ is a completely bounded 
projection from $\tilde Y$ into $\tilde X$. 
Thus $X$  has the CSCP.

(c) $\To$ (a): 
The proof is the same. 
Indeed, letting $\tilde X$, $\tilde Y$ and $\tilde T$ be as above, $Y$ 
is also separable and locally reflexive, hence since $X$ is thus locally 
complemented in $Y$ by assumption, $X$ is completely complemented in $Y$, 
so $\tilde X$ is completely complemented in $\tilde Y$ as shown above. 
\enddemo

(b) $\To$ (c) is an immediate consequence of the following:  

\proclaim{Lemma A} 
Let $Y$ and $Z$ be locally reflexive operator spaces and $X$ be a locally 
complemented subspace of $Y$. 
Let $T:X\to Z$ be a completely bounded linear operator. 
Then there exists a locally reflexive operator space $W\supset Z$ and a 
completely bounded linear operator $\tilde T:Y\to W$ with $\tilde T|X=T$. 
\endproclaim

\demo{Remark}
We obtain that in fact $\tilde T$ may be chosen with $\|\tilde T\|_{\cb} 
= \|T\|_{\cb}$ and $W$ is $\beta$-locally reflexive, where if $X$ is 
$(C+\ep)$-locally complemented in $Y$ for all $\ep>0$, $Y$ is 
$\lambda$-locally reflexive, and $Z$ is $\gamma$ locally reflexive, then 
$$\beta = (\max \{\gamma,\lambda\}) (C+1)\ .$$
\enddemo

We first deduce (b) $\To$ (c) of Theorem~1.4 from the Lemma. 
Simply choose $W\supset X$ locally reflexive and $T':\tilde Y\to W$ a 
completely bounded map extending $T$, by the Lemma, and then let 
$P:W\to X$ be a completely bounded surjective projection; 
$\tilde T= PT'$ is then the desired extension of $T$.\qed

\demo{Comment} 
Of course the implication (b) $\To$ (a) is superfluous in this entire 
argument, but its proof is considerably simpler than that of Lemma~A.
\enddemo 

To prove Lemma A, we require 

\proclaim{Lemma B} 
Let $X$ be a locally complemented subspace of a locally reflexive operator 
space $Y$. 
Then $Y/X$ is locally reflexive.
\endproclaim

\demo{Proof} 
Assume $X$ is $(C+\ep)$-locally complemented in $Y$, for all $\ep>0$. 
By the Remark following the statement of Theorem~1.4, choose a completely 
bounded projection $P:Y^{**}\to X^{**}$ with $\|P\|_{\cb} \le C$; 
let $E$ denote the null space of $P$. 
It follows that if $\Pi :Y\to Y/X$ is the quotient map, then $\Pi^{**}|E
\to (Y/X)^{**}$ is a complete surjective isomorphism with 
$$\|(\Pi^{**}|E)^{-1}\|_{\cb} \le C+1 \ .$$
Now let $G$ be a finite-dimensional subspace of $(Y/X)^{**}$ and $F$ be a 
finite dimensional subspace of $(Y/X)^* = X^\bot$; 
set $\tilde G = (\Pi^{**}|E)^{-1}G$. 
Now assuming that $Y$ is $\lambda$-locally reflexive, given $\ep>0$, choose 
$\tilde T:\tilde G\to Y$ with 
$$\|\tilde T\|_{\cb} < \lambda+\ep\text{ and } 
\langle \tilde T\tilde g,f\rangle = \langle \tilde g,f\rangle
\text{ for all } \tilde g\in\tilde G,\ f\in F\ .$$
Finally, define $T:G\to Y/X$ by 
$$T= \Pi \tilde T(\Pi^{**}|E)^{-1}\ .$$
Evidently then 
$$\|T\|_{\cb} \le (\lambda+\ep)(C+1)\ .$$

Finally, suppose $g\in G$ and $f\in F$. 
Then 
$$\eqalignno{
\langle \Pi \tilde T(\Pi^{**}|E)^{-1} g,f\rangle
&=\langle\tilde T(\Pi^{**}|E)^{-1}g,f\rangle\ \text{ (since $f\in X^\bot$)}\cr
&= \langle (\Pi^{**}|E)^{-1} g,f\rangle\cr
&= \langle \Pi^{**} (\Pi^{**}|E)^{-1}g,f\rangle\ \text{ (again since 
$f\in X^\bot$)}\cr
&= \langle g,f\rangle&\qed\cr}$$
\enddemo

\demo{Remarks} 
1. Of course we obtain that $Y/X$ is $\lambda (C+1)$-locally reflexive. 
Actually, if we assume instead that  $X^{**}$ is $C$-completely 
cocomplemented in $Y^{**}$, we have that $Y/X$ is $(\lambda C)$-locally 
reflexive. 
In particular, {\it if $Y$ is $1$-locally reflexive and $X^{**}$ is 
completely contractively cocomplemented in $Y^{**}$, $Y/X$ is $1$-locally 
reflexive.}

2. After the ``final'' draft of our paper was finished, we learned 
that S-H.~Kye and Z-J.~Ruan had already obtained a variant of Lemma~B in 
\cite{KR} (see Proposition~5.4 there), as well as an interesting converse. 
The work in \cite{KR} contains moreover some remarkable characterizations 
of $\lambda$-injectivity for dual operator spaces. 
\enddemo

\demo{Proof of Lemma A} 
We use the construction of G.~Pisier mentioned in the Remarks following 
the statement of Proposition~1.3. 
Let $X,Y,Z$ and $T$ be as in Proposition~1.3; we may assume that 
$\|T\|_{\cb}=1$. 
Assume then that $C$, $\lambda$ and $\gamma$ are as in the Remark following 
the statement of Lemma~A. 
Let $E= (Y\oplus Z)_\infty$ and $\Gamma = \{x\oplus -Tx:x\in X\}\subset E$; 
let $W = E/\Gamma$, $\Pi : E\to W$ the quotient map, and define 
$j:Z\to W$ and $\tilde T:Y\to W$ by $j(z) = \pi (0\oplus z)$ and 
$\tilde T(y) = \pi (y\oplus 0)$ for all $z\in Z$ and $y\in Y$. 
Then $j$ is a complete (into) isometry and $jT= \tilde Ti$, where 
$i:X\to Y$ is the inclusion map. 
Thus $\tilde T$ is the desired extension with $\|\tilde T\|_{\cb}=1$ also. 

Now we have that $E$ is $\max \{\lambda,\gamma\}$ locally reflexive. 
We claim that $\Gamma$ {\it is locally complemented in $E$}. 
To see this, let $P$ be the  natural projection of $E$ onto $Y$ 
with nullspace $Z$ and let $\Lambda$ be a linear subspace of $E$ containing 
$\Gamma$ with $\Gamma$ of finite-codimension in $\Lambda$. 
Now $P\Gamma=X$ and so $X$ is of finite-codimension in $P\Lambda$. 
Thus given $\ep>0$, there is a surjective linear projection 
$Q:P\Lambda\to X$ with $\|Q\|_{\cb} <C+\ep$. 
Now defining $U:X\to \Gamma$ by $U(x) = x\oplus Tx$ for all $x\in X$, 
also $\|U\|_{\cb} =1$ and of course $UX=\Gamma$; we claim that 
$R\defeq UQP|\Lambda$ is the desired projection onto $\Gamma$. 
Clearly $RE\subset\Gamma$. 
But let $x\oplus -Tx$ be a a typical element of $\Gamma$. 
Then 
$$(UQP) (x\oplus -Tx) = UQx= Ux = x\oplus -Tx\ .$$
We have thus proved that also $\Gamma$ is $(C+1)$-locally complemented. 
Thus $W$ is indeed locally reflexive by Lemma~B.\qed
\enddemo

We now proceed with the proof of Theorem 1.1. 
For the remainder of this section, we assume $H$ is separable; 
we identify $H$ with $\ell^2$ and $B(H)$ with $\M_\infty$, the set of all 
infinite matrices yielding bounded linear operators on $H$. 
$\N \subset \M_\infty$ will be called  a {\it special copy\/} of $B(H)$ 
provided there exist infinite pairwise-disjoint subsets $M_1,M_2,\ldots$ 
of $N$ so that letting $M = \bigcup_{j=1}^\infty M_j$, 
and letting $m_1^j,m_2^j,\ldots$ be an increasing enumeration of 
$M_j$ for each $j$,
then $A= (a_{ij})$ in $\M_\infty$ belongs to $\N$ provided 
\roster
\item"(i)" $a_{ij} = 0$ if $i$ or $j\notin M$
\item"(ii)" for all $i,j\in M$, if $i\in M_r$ and $j\in M_s$, 
there exist numbers $b_{rs}$ so that 
$a_{ij} = b_{rs}$ if $i=m_k^r$ and $j=m_k^s$ for some $k$, 
$a_{ij}=0$ if $i=m_k^r$ and $j=m_\ell^s$ with $k\ne \ell$. 
\endroster
Evidently $\N$ is then a WOT-closed $*$-subalgebra of $B(H)$, 
$*$-isomorphic to $B(H)$. 
(We could also define $\N$ ``intrinsically'' as follows: 
let $e_1,e_2,\ldots$ be the usual basis for $\ell^2$; let 
$H_j = [e_i]_{i\in M_j}$. 
For each $i,j$ let $E_{ij}$ be the partial isometry in $B(H)$ with initial 
domain $H_j$ and final domain $H_i$, so that 
$E_{ij}(e_{m_\ell^j}) = e_{m_\ell^i}$, $\ell=1,2,\ldots$. 
Then $\N$ is the WOT closed linear span of the $E_{ij}$'s.)

\proclaim{Lemma 1.5} 
Let $Z$ be a separable closed subspace of $B(H)$. 
There exists a special copy $\N$ of $B(H)$ so that $Z\oplus \N$ is a complete 
direct decomposition.
\endproclaim

\remark{Remark} 
In fact we show that letting $Q$ be the projection from $Z\oplus \N$ 
onto $\N$ with kernel $Z$, then 
$$\|Q\|_{\cb} \le 2\ .
\tag 14$$
\endremark

Let us first give the 

\demo{Proof of Theorem 1.1}
Let $X,Y$ and $T$ be as in the statement of 1.1. 
Since $B(H)$ is 1-injective, we may 
choose $T':Y\to B(H)$ a completely bounded linear map extending $T$ 
with $\|T'\|_{\cb} = \|T\|_{\cb}$. 
Now set $Z= \overline{T'Y}$. 
Choose $\N$ a special copy of $B(H)$ satisfying (14) (by Lemma~1.5). 
Since $\N$ is completely isometric to $B(H)$ and $Y/X$ is separable, we 
may choose $V:Y/X\to \N$ a complete (into) isometry. 
Let $\Pi :Y\to Y/X$ denote the quotient map and define $\tilde T:Y\to B(H)$ by 
$$\tilde T y = T'y + V\Pi y\ \text{ for all }\ y\in Y\ .
\tag 15$$
Now letting $W= Y/X$, we have that the map $U:Z\oplus W\to Z+V(W)$ is 
a complete isomorphism, where 
$$U(z\oplus w) = z+V(w) \ \text{ for all }\ z\in Z\ \text{ and }\ w\in W\ .$$
Indeed, letting $I$ denote the identity on $\bK$ and $P$ the projection from 
$Z\oplus W$ onto $W$ with kernel $Z$ and $R= Id-P$, then $U= VP+R$, so 
$$\|I\otimes U\| \le \|I\otimes VP\| + \| I\otimes R\| \le 2\ .
\tag 16$$
Also if $\oQ$ is the projection from $Z+V(W)$ onto $V(W)$ with kernel $Z$, 
and $\tau \in \bK \otimes V(W)$, then 
$$\align 
\|I\otimes U^{-1}(\tau)\| 
& = \|I\otimes U^{-1} \oQ\tau + I\otimes U^{-1} (I-\oQ)\tau\| 
\tag 17\\
& = \max \left\{\|I\otimes U^{-1} \oQ(\tau)\|, 
\|I\otimes U^{-1}(I-\oQ)\tau\|\right\} \\
&\le 3\|\tau\|\ \text{ by (14).}
\endalign$$
Thus we have by (16) and (17) 
$$\|U\|_{\cb} \|U^{-1}\|_{\cb} \le 6\ .
\tag 18$$
Now if we instead define $\skew4\tilde{\tilde T} : Y\to Z\oplus Y/Z$ by 
$\skew4\tilde{\tilde T} (y) = T'y\oplus \Pi y$, then by the 
proof of Proposition~1.3, $\skew4\tilde{\tilde T}$ is a complete into 
isomorphism, hence also $\tilde T = U\skew4\tilde{\tilde T}$ is a complete 
into isomorphism, and we have by (8) that $(*)$ holds.\qed
\enddemo

We now proceed with the proof of Lemma 1.5. 
For $M\subset \nat$, let $H_M =  [e_i:i\in M]$ and let 
$B_M = PB(H)P$ where $P$ is the orthogonal projection on $H_M$. 

\proclaim{Lemma 1.6} 
Let $G$ be a separable subspace of $\bK^\bot$. 
There exists an infinite $M$ with $G\perp B_M$.
\endproclaim

\demo{Proof} 
Say that $L,M$ are almost disjoint if $(L\sim M)\cup (M\sim L)$ is finite. 
Now if $L,M$ are almost disjoint subsets of $\nat$ and $\mu\in \bK^\bot$, then 
$$\|\mu\| \ge \|\mu |B_L\| + \|\mu|B_M\| \ .
\tag 19$$ 
Indeed, we may choose $L'\subset  L$, $M'\subset M$ with $L',M'$ disjoint 
and $L\sim L'$, $M\sim M'$ finite. 
Since $\mu\in \bK^\bot$, we have e.g., $\mu|B_{L\sim L'} =0$, so 
$\mu |H_L = \mu |B_{L'}$, $\mu |B_M = \mu |B_{M'}$, whence 
$$\|\mu |B_L\| + \|\mu |B_M\| 
= \|\mu |B_{L'} \| + \|\mu |B_{M'}\| 
= \|\mu |B_{L'\cup M'} \| \le \|\mu \| .
\tag 20$$

It follows immediately that if $M_1,\ldots, M_n$ are pairwise almost 
disjoint subsets of $\nat$, then 
$$\|\mu \| \ge \sum_{i=1}^n \|\mu (B_{n_i})\|\ .
\tag 21$$
Now by an ancient classical fact, there exists an uncountable family 
$(M_\alpha)_{\alpha\in\Gamma}$ of almost disjoint infinite subsets of $\nat$. 
(21) yields that for $\mu \in \bK^\bot$, 
$$\mu |B_{M_\alpha} =0\ \text{ for all but countably many $\alpha$'s.}
\tag 22$$ 
Since $G$ is separable, it now follows that for some $\alpha$, 
$M|B_{M_\alpha} = 0$ for all $M\in G$, yielding 1.6.\qed
\enddemo

We need one more ingredient to obtain Lemma 1.5. 
Let $\Ca$ denote the Calkin algebra $B(H)/\bK$, and now let 
$\Pi :B(H) \to \Ca$ denote the quotient map. 

\proclaim{Lemma 1.7} 
Let $\N$ be a special copy of $B(H)$. 
Then $\Pi |\N$ is an (into) $*$-isomorphism. 
In particular, $\Pi|\N$ is a complete isometry.
\endproclaim

\demo{Proof} 
Since $\Pi$ is a $*$-isomorphism, we need only observe that $\N$ contains 
no non-zero compact operators, which is trivial since $E_{ij}$ is a 
non-finite rank partial isometry for all $i,j$ as in the alternate 
definition of special $\N$.\qed
\enddemo

We are finally prepared for the 

\demo{Proof of Lemma 1.5} 
Since $Z$ is separable, so also is $\bK \otimes \overline{\Pi Z}$; thus  we 
may choose $G$ a separable subspace of $\bK^\bot$ so that 
$$\bK \otimes G \text{ isometrically norms } \bK\otimes \overline{\Pi Z} 
\text{ via the canonical pairing.}
\tag 23$$
Now by Lemma 1.6, choose $M$ an infinite subset of $\nat$ so that 
$$G\perp B_M\ .
\tag 24$$ 
Finally, let $M_1,M_2,\ldots$ be infinite disjoint sets with $M=\bigcup M_j$, 
and let $\N$ be the special copy of $B(H)$ corresponding to the $M_j$'s. 
Of course then  $\N \subset B_M$, and so by (24), 
$$g(x) = g(\Pi x)=0\ \text{ for all } g\in G\text{ and } x\in \N\ .
\tag 25$$ 
Now let $\tau$ be an element of $\bK\otimes (Z\oplus \N)$, say 
$$\tau = \sum_{i=1}^k L_i \otimes (z_i \oplus x_i)\ \text{ where }\ 
L_i\in \bK\ ,\ z_i\in Z\text{ and } x_i\in\N \text{ for all } i\ .$$
It follows from Lemma 1.7 that then 
$$\|\sum L_i\otimes x_i\| = \|\sum L_i \otimes \Pi x_i\|\ .
\tag 26$$
Now by (23), given $\ep>0$, we may choose $S_1,\ldots,S_\ell$ in $\bK$ 
and $g_1,\ldots,g_\ell$ in $G$ so that  $\|\sum S_j \otimes g_j\|=1$ and 
$$\align
\|\sum L_i \otimes \Pi (z_i) \| 
& \le (1+\ep)|\langle \sum S_j\otimes g_j,\sum L_i\otimes \Pi z_i\rangle | \\
&= (1+\ep)|\langle \sum S_j\otimes g_j,\sum_i L_i\otimes 
\Pi (z_i\oplus x_i)\rangle |  \text{ (by  (25)}
\tag 26\\
&\le (1+\ep) \|\sum L_i \otimes \Pi (z_i\oplus x_i)\| 
\endalign$$
(where $\langle \sum S_j\otimes g_j,\sum L_i\otimes T_i\rangle 
= \sum_{i,j} g_i (T_i) S_j\otimes L_i$). 

\noindent Thus 
$$\align
\|\sum L_i\otimes x_i\| 
& = \|\sum L_i\otimes\Pi x_i\|\text{ (by Lemma~1.7)}
\tag 27\\
&\le (2+\ep) \|\sum L_i\otimes\Pi (z_i\oplus x_i)\|\text{ by (26)} \\
&\le (2+\ep) \|\sum L_i\otimes (z_i\oplus x_i)\|\ .
\endalign$$
Since $\ep>0$ is arbitrary, we have indeed proved that if $Q(z+x)=x$ for all 
$z\in Z$ and $x\in\N$, then $\|Q\|_{\cb} \le 2$.\qed
\enddemo

We conclude this section with a considerable strengthening of Corollary~1.4. 

\proclaim{Corollary 1.8} 
Let $X$ be a locally reflexive separable operator space so that $X^{**}$ 
is completely isomorphic to $B(H)$. 
Then $X$ has the {\rm CSCP} provided $X$ is completely complemented in every 
separable locally reflexive $Y$ with 
$$X\subset Y \subset X^{**}\ .$$
\endproclaim

\demo{Proof} 
Let $Z\subset W$ be operator spaces with $W$ separable locally reflexive 
and $Z$ completely isomorphic to $X$. 
We must show that $Z$ is completely complemented in $W$. 
Let $T:Z\to X$ be a complete surjective isomorphism. 
By Theorem~1.1, since $X^{**}$ is completely isomorphic to $B(H)$, we 
may choose $Y$, $X\subset Y\subset X^{**}$ and $\tilde T:W\to Y$ a complete 
surjective isomorphism extending $T$. 
Then $Y$ is locally reflexive, hence there is a completely bounded projection 
$P$ from $Y$ onto $X$. 
Then $Q= T^{-1}P \tilde T$ is a completely bounded projection from $W$ onto 
$Z$.\qed
\enddemo 

\head {\smc Section} 2\\
{\bf An operator space construction}\\
{\bf on certain subspaces of $\M_\infty$}\endhead

\definition{Definition} 
Let $W\subset\nat \times\nat$. 
$M_W$ denotes all $A$ in $\M_\infty$ with $a_{ij} = 0$ if $(i,j)\notin W$. 
$\bK_W$ denotes $M_W\cap \bK$. 
We {\it define\/} an operation on all $\nat\times\nat$ matrices, denoted $P_W$, 
as follows: 
for any $A$, $P_W A=B$ where $b_{ij}= a_{ij}$ if $(i,j) \in W$, $b_{ij}=0$ 
otherwise. 
In case $P_W\M_\infty \subset \M_\infty$, it follows immediately that 
$P_W|\M_\infty$ is bounded; then we set 
$$\|P_W\| = \|P_W|\M_\infty\|\ \text{ and }\  
\|P_W\|_{\cb} = \|P_W|\M_\infty\|_{\cb}$$
(which a-priori might be infinite). 
\enddefinition

We may now formulate the main result of this section, which involves 
the construction of a certain $W$ for which $P_W$ is completely bounded. 

\proclaim{Theorem 2.1} 
Let $Y$ be a separable closed subspace of $\M_\infty$ with $\bK \subset Y$, 
and let $\ep>0$. 
There exists an absolute constant $C$, a subset $W$ of $\nat\times\nat$, and 
a subspace $\tilde Y$ of $B(H)$ satisfying the following: 
\roster
\item"(i)" $\|P_W\|_{\cb} \le 2$.
\item"(ii)" $d_{\cb} (\bK_W, \bK_0) \le C$.
\item"(iii)" There is a complete surjective isomorphism $T:Y\to\tilde Y$
with 
\vskip1pt
\itemitem{\rm (a)} $\|T\|_{\cb} \|T^{-1}\|_{\cb} \le 1+\ep$
\itemitem{\rm (b)} $\|Ty-y\| \le \ep \|y\|$ for all $y\in Y$ 
\itemitem{\rm (c)} $T|K = I|K$.
\vskip1pt
\item"(iv)" $\tilde Y$ is invariant under $P_W$.
\item"(v)" $P_{\sim W} \tilde Y\subset \bK$.
\endroster
\endproclaim

Before proving this result, we give two applications.

\proclaim{Theorem 2.2} 
$\bK$ has the CSCP.
\endproclaim

\demo{Proof} 
By Corollary 1.8, it suffices to prove that if $Y$ is locally reflexive 
separable and 
$$\bK \subset Y \subset \M_\infty\ ,
\tag 28$$ 
then $\bK$ is completely complemented in $Y$. 
Now let $\ep>0$, and choose $W$ and $\tilde Y$ as in Theorem~2.1; also let 
$T$ satisfy (iii) of Theorem~2.1. 
Then by (iii)(c) and (iv), 
$$\bK\subset \tilde Y\ \text{ and }\ \bK_W\subset \tilde Y\ .
\tag 29$$ 
By (ii) and (iii), $\bK_W$ is completely isomorphic to $K_0$ and $\tilde Y$ 
is separable locally reflexive. 
Hence by the results in \cite{R}, there is a completely bounded projection 
$Q$ from $\tilde Y$ onto $\bK_W$. 
Then using (29) and (iv), (v) of Theorem~2.1, $\tilde P$ is a completely 
bounded projection from $\tilde Y$ onto $\bK$, where 
$$\tilde P = (QP_W + P_{\sim W})|\tilde Y\ .
\tag 30$$
Finally, $P\defeq T^{-1} \tilde P T$ is a completely bounded projection 
from $Y$ onto $\bK$, completing the proof.\qed
\enddemo

For our second application, we briefly introduce the concept to be 
developed in the next section. 

\definition{Definition}
An operator space $Z$ has the {\it Mixed Separable Extension Property\/} (MSEP) 
if for all separable operator spaces $Y$, subspaces $X$, and completely 
bounded maps $T:X\to Z$, there exists a bounded linear map $\tilde T:Y\to Z$ 
extending $T$.
\enddefinition

As we show in Section 3, a separable $Z$ has the MSEP iff $Z$ is complemented 
in $Y$ for every separable operator space $Y$ with $Z\subset Y\subset B(H)$ 
(for $Z\hookrightarrow B(H)$ a fixed complete embedding). 
As noted in the introduction, we do not know if $\bK$ has this property.
The next result reduces this problem to $\bK_0$.

\proclaim{Proposition 2.3} 
$\bK$ has the {\rm MSEP} if and only if $\bK_0$ does.
\endproclaim

\demo{Proof} 
If $\bK$ has the MSEP so does $\bK_0$, because it is completely complemented 
in $\bK$. 
For the non-trivial implication, suppose $\bK_0$ has the MSEP and let $Y$ 
be a separable operator space satisfying (28). 
Again, for fixed $\ep>0$, choose $W,\tilde Y$, and $T$ as in Theorem~2.1; 
then choose $Q$ a bounded linear projection from $\tilde Y$ onto $\bK_W$. 
Again, letting $\tilde P$ and $P$ be as in the proof of 2.2, it now follows 
that $P$ is a bounded linear projection from $Y$ onto $\bK$.\qed
\enddemo

We now proceed with the proof of Theorem 2.1. 
We first isolate part of the proof in the following result:  

\proclaim{Lemma 2.4} 
Let $W\subset \nat\times\nat$ be described as follows: 
there exists $(m_j)$ in $\nat \cup \{0\}$ that $1=m_0<m_1<m_2<\cdots$, 
$m_{j+1} -m_j\to\infty$, so that for all $(i,j)\in \nat\times\nat$, 
$(i,j) \in W$ iff the following are all satisfied for some $k=1,2,\ldots$;
\roster
\item"(a)" $m_{k-1} <i\le m_k$ and $j\le m_{k+1}$;
\item"(b)" $m_{k-1} < j\le m_k$ and $i\le m_{k+1}$; 
\item"(c)" $i=1$ and $j\le m_1$ or $j=1$ and $i\le m_1$.
\endroster
Then $\|P_W\|_{\cb} \le2$ and $d_{\cb}(\bK_W,\bK_0)\le C$ for some 
absolute constant $C$.
\endproclaim

In the following we use interval notation to denote intervals in 
$\nat\cup\{0\}$.

\demo{Proof}
Let $A_j = (m_{2j-1},m_{2j+1}]\times (m_{2j-1},m_{2j+1}]$ for $j\ge 1$,\newline
\hphantom{{\it Proof.} Let }$A_0 = [m_0,m_1]\times [m_0,m_1]$.\newline
\hphantom{\it Proof.} Let $B_j = (m_{2j},m_{2(j+1)}] \times 
(m_{2j},m_{2(j+1)}]$ for $j\ge0$.\newline
\hphantom{\it Proof.} Let $\ A = \bigcup\limits_{j=0}^\infty A_j$, \quad 
$B= \bigcup\limits_{j=0}^\infty B_j$. 

We claim that 
$$W= A\cup B\ .
\tag 31$$
The following diagram intuitively illustrates why this is so: 
the heavy lines denote the $A_j$'s, the dotted lines denote the 
$B_j$'s. $\bigcup_{j=0}^2 B_j \sim A$ is shaded in the diagram; 
the regularity of $\bigcup_{j=0}^\infty B_j\sim A$ used in showing 
that $P_W$ is completely bounded. 
$$\epsfysize=2.3truein\epsfbox{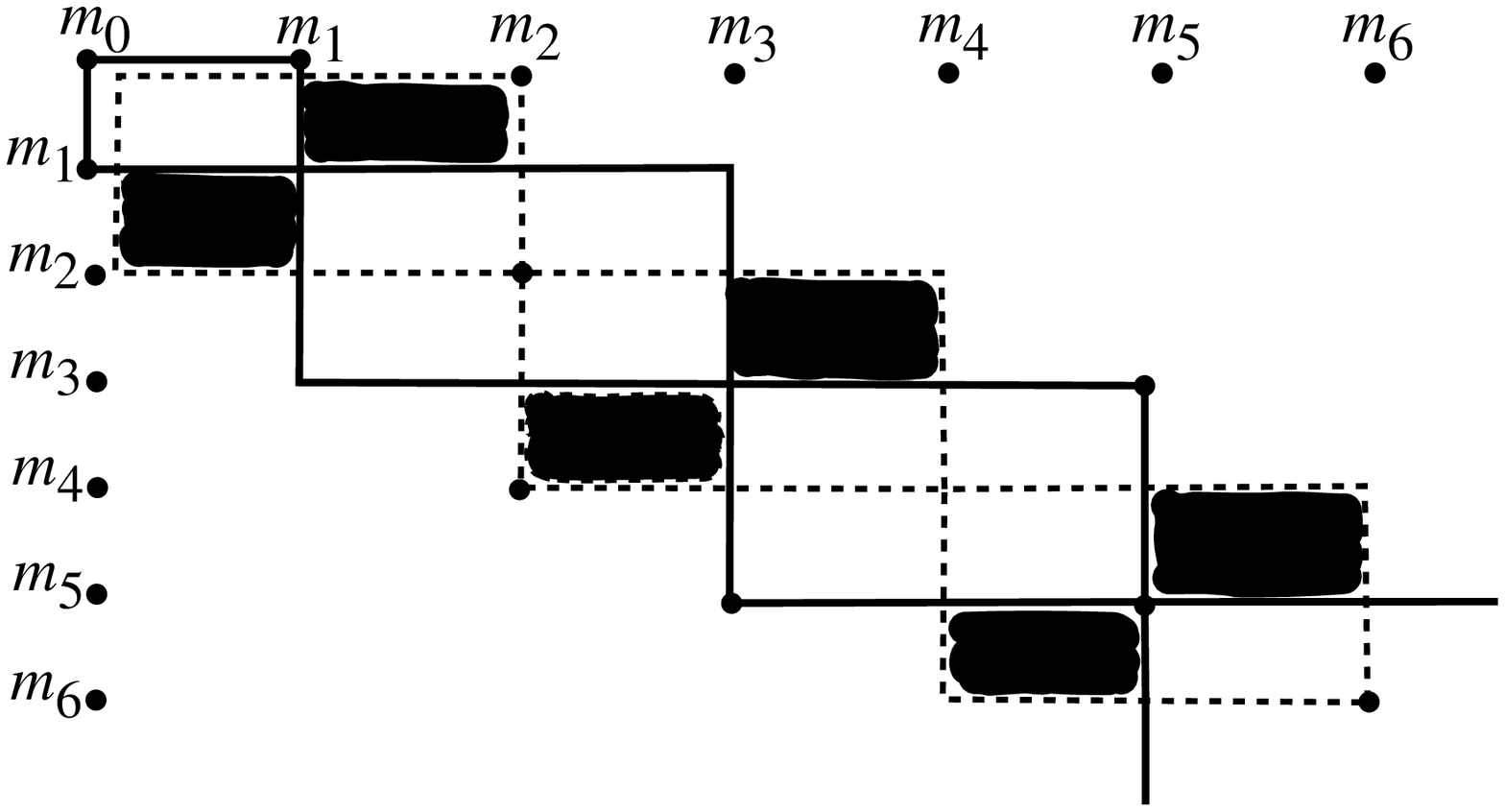}$$
First, if $(i,j)$ satisfy (c) of 2.4, then $(i,j)\in A_0$. 
Now suppose $(i,j)\in W$, and $i>1$, $j>1$. 
Now if $i\le m_1$ or $j\le m_1$, then
$(i,j)\in (m_0,m_2)\times (m_0,m_2]\subset B$ by (a) and (b) of 2.4. 

Suppose then $i>m_1$ and $j>m_1$. 
But then if $i\le m_2$ or $j\le m_2$, $(i,j)\in (m_1,m_3]\times [m_1,m_3]
\subset A$. 
Continuing by induction, we obtain that 
$$W\subset A\cup B\ .
\tag 32$$ 
Next suppose $(i,j)\in [m_0,m_j]\times [m_0,m_1]$. 
By (c), we may assume $i>1$ and $j>1$. 
But then $i$ and $j$ satisfy both (a) and (b) for $k=1$, so $(i,j)\in W$. 
Now suppose $(i,j) \in (m_0,m_2]\times [m_0,j_2]\sim [m_0,m_1]\times 
[m_0,m_1]$. 
Then if $1<i\le m_1$, $m_1<j\le m_2$, whence (a) holds for $k=1$ and also 
(b) holds vacuously for $k=1$, while (a)  holds vacuously for $k=2$ and 
(b) holds for $k=2$.

If $m_1<i\le m_2$, we get that $1 <j\le m_1$, so by symmetry again (a) and 
(b) both hold for $k=1$ and $k=2$. 
Thus $(m_0,m_2] \times (m_0,m_2]\subset W$. 
Carrying this one more step for the pattern, now suppose 
$$(i,j) \in (m_1,m_3] \times (m_1,m_3] \sim (m_0,m_2] \times (m_0,m_2]\ .$$
Thus if $m_1 < i\le m_2$, $m_2<j\le m_3$, whence (a) holds for $k=2$, 
vacuously for $k=3$, and (b) holds for $k=2$ and $k=3$.

If $m_2 <i\le m_3$, then $m_1 <j\le m_2$, so again (a) and (b) hold 
for $k=2$ and $k=3$ by symmetry. 

Thus by induction, we obtain that 
$$(m_{2j},m_{2j+2}] \times (m_{2j},m_{2j+2}]\ \text{ and }\ 
(m_{2j+1},m_{2j+3}] \times (m_{2j+1},m_{2j+3}]\subset W$$ 
for all $j$, whence 
$$A\cup B\subset W\ .
\tag 33$$
Of course (31) is now established via (32) and (33). 

Now it is evident that $\|P_E\|_{\cb} =1$ for $E= A,B$, and $A\cap B$. 
This gives the ``easy'' estimate 
$$\|P_W\|_{\cb} \le 3\ ,
\tag 34$$ 
since $P_W = P_A + P_B - P_{A\cap B}$. 
We are indebted to T.~Schlumprecht for the following better estimate:
$$\|P_{B\sim A}\|_{\cb} =1\ .
\tag 35$$ 
To see this, fix $j\ge 0$. 
Then 
$$\align
&(m_{2j},m_{2j+2}] \times (m_{2j},m_{2j+2}]\sim A \tag 36\\
&\qquad = (m_{2j},m_{2j+1}] \times (m_{2j+1},m_{2j+2}] 
\cup (m_{2j+1},m_{2j+2}] \times (m_{2j}, m_{2j+1}]\ . 
\endalign$$
Resorting to a simple picture, we thus have that he matrices in 
$M_{B_j\sim A}$ have the form 
$$T = \left[ \matrix 0&C\\ \noalign{\vskip6pt} D&0\endmatrix\right]$$
whence $\|T\| = \max \{ \|C\|,\|D\|\}$ and so 
$\|P_{B_j\sim A}\|_{\cb} =1$. 
Since  $B$ is the union of the disjoint blocks $B_1,B_2,\ldots,$ (35) follows. 

It remains to prove the final assertion of 2.4. 
In the following, the letter ``$c$'' denotes absolute constants, which 
may vary from line to line.

First, via the Pe{\l}czy\'nski decomposition method, we obtain the following

\example{Fact} 
Let $(n_j)$ be a sequence of positive integers with $\sup_j n_j=\infty$. 
Then 
$$d_{\cb}  \left( \Big( \bigoplus_{j=1}^\infty M_{n_j} \Big)_{c_0}, 
\bK_0\right) \le c\ .$$
(In fact, here one may take $c=2$.) 
It then follows immediately that 
$$d_{\cb} (\bK_A ,\bK_0) \le c
\tag 37$$
(for $c$ in the Fact).  
Next, we define $C,D$ by 
$$\align
C& = \bigcup_{j=0}^\infty (m_{2j},m_{2j+1}] \times (m_{2j+1},m_{2j+2}]
\tag 38i \\
D& = \bigcup_{j=0}^\infty (m_{2j+1}, m_{2j+2}] \times (m_{2j},m_{2j+1}]\ .
\tag 38ii
\endalign$$
Again by the Fact, we obtain 
$$d_{\cb} (\bK_E,\bK_0) \le c\ \text{ for }\ E= C\text{ or }D\ .
\tag 39$$
Finally, we have 
$$K_W = K_A \oplus K_C \oplus K_D\ .
\tag 40$$
Indeed, $\|P_E\|_{\cb} =1$ for $E=A$, $C$ or $D$, and 
$W= A\cup C\cup D$ by (31) and (38). 
We then have that 
$$d_{\cb} (\bK_W, (\bK_A \oplus \bK_C \oplus \bK_D)_\infty) \le 3
\tag 41$$
and by the Fact 
$$d_{\cb} ((\bK_A \oplus \bK_C \oplus \bK_D)_\infty, \bK_0) \le c
\tag 42$$ 
completing the proof.\qed
\endexample
\enddemo

\remark{Remark} 
The following intriguing problem arises: 
{\it characterize the sets $W\subset \nat\times\nat$ so that $P_W$ 
is bounded\/}. 
A related problem:
{\it if $P_W$ is bounded, is it completely bounded?}
\endremark

\demo{Proof of Theorem 2.1}

Let $\Pi :Y\to Y/\bK$ be the quotient map; without loss of generality, 
$Y/\bK$ is infinite-dimensional. 
Choose $y_1,y_2,\ldots$ in $Y$ so that $(y_j)$ is bounded and $(\Pi y_j)$ 
is a bounded biorthogonal system spanning $\Pi (Y)$. 
Thus, we may choose $M<\infty$ and $(y_j^*)$ in $Y^\bot$ so that for 
all $j$ and k, 
$$y_j^* (y_k) = \delta_{jk}\quad ,\quad 
\|y_j^*\| \le M\text{ and } \|y_j\|\le M\ .
\tag 43$$
Let $\ep>0$ be given and set $m_0=0$. 
We shall construct a sequence $(m_j)$ in $\nat$ and certain sequences 
$(y_i^{(j)})$ in $Y$. 
\enddemo

\demo{Step 1} 
Choose $m_1\in\nat$ so that 
$$\|P_{\{1\} \times (m_1,\infty)} y_1\| 
+ \|P_{(m_1,\infty)\times \{1\}} y_1\| < \frac{\ep}2\ .
\tag 44$$
Now define $y_1^{(1)} = y_1$ and $y_j^{(1)} = P_{(m_1,\infty)\times 
(m_1,\infty)} y_j$ for all $j>1$.
\enddemo

\demo{Step 2} 
Choose $m_2 >m_1$ so that 
$$\| P_{[1,m_1]\times (m_2,\infty)} y_i^{(1)}\| 
+ \|P_{(m_2,\infty)\times [1,m_1]} y_i^{(1)}\| < \frac{\ep}{2^2} 
\tag 45$$ 
for $1\le i\le 2$. 

Now set $y_i^{(1)} = y_i^{(2)}$ for $i\le 2$, 
$$y_i^{(2)} = P_{(m_2,\infty)\times (m_2,\infty)} y_i^{(1)})\ \text{ for } 
i>2\ .$$
\enddemo

\demo{Step $j$}
Suppose $j>2$ and $m_1<\cdots < m_{j-1}$, $(y_i^{(s)})_{i=1}^\infty$ have 
been chosen, for all $1\le s\le j-1$. 
Choose $m_j>m_{j-1}$ so that 
$$\|P_{[1,m_{j-1}]\times (m_j,\infty)} y_i^{(j-1)}\| 
+ \| P_{(m_j,\infty)\times [1,m_{j-1}]} y_i^{j-1}\| < \frac{\ep}2\ ,
\tag 46$$ 
for all $i$, $1\le i\le j$. 

Now set $y_i^{(j)} = y_i^{(j-1)}$ for $1\le i\le j$, then set 
$$y_i^{(j)} = P_{(m_j,\infty)\times (m_j,\infty)} y_i^{(j-1)}
\ \text{ for all } i>j\ .
\tag 47$$
This completes the inductive construction. 
Then we have that for all $j\ge1$, 
$$y_i^{(j)} = P_{(m_j,\infty)\times (m_j,\infty)} y_i\ \text{ for all }\ i>j\ .
\tag 48$$
Hence 
$$y_j^{(j)} = P_{(m_{j-1},\infty)\times (m_{j-1},\infty)} y_j
\ \text{ for all }\ j>1 \text{ and } y_1^{(1)} =y_1\ .
\tag 49$$
Thus it follows that $y_j-y_j^{(j)}$ is a finite rank operator for all $j$, 
whence 
$$\Pi y_j^{(j)}=\Pi y_j\ \text{ and }\ \|y_j^{(j)}\|\le M\text{ for all }j\ .
\tag 50$$

Now let $W$ be defined as  in Lemma 2.4. 
That is, instead defining $C_0,C_1,C_2,\ldots$ 
$$\align
C_k &= (m_{k-1},m_k]\times (m_{k+1},\infty) 
\cup (m_{k+1},\infty)\times (m_{k-1},m_k]\ \text{ for } k\ge1\ ,\tag 51\\
C_0 & = \{1\} \times (m_1,\infty) \cup (m_1,\infty) \times\{1\}\ ,\tag 52
\endalign$$
then 
$$W= \sim C \ \text{ where }\ C= \bigcup_{j=0}^\infty C_j\ .
\tag 53$$
Now fix $1\le i$. 
Then by (44)--(46), (49) and (51)--(52), 
$$\|P_{C_{j-1}} y_i^{(i)}\| < \frac{\ep}{2^j} \ \text{ if }\ j\ge i\ .
\tag 54$$
But by (49), 
$$P_{C_{j-1}} y_i^i = 0\ \text{ if }\ i>j\ .
\tag 55$$
We thus obtain for all $i$, from (54) and (55), that 
$$\|P_C y_i^{(i)}\| = \Big\| \sum_{j=i}^\infty P_{C_{j-1}} y_i^{(i)}
\Big\| \le \sum_{j=i}^\infty \|P_{C_{j-1}} y_i^i\| 
< \frac{\ep}{2^{i-1}}\ .
\tag 56$$

We next define the operator $T$ specified in the statement of 2.1. 
First, let $Y_0$ denote the {\it linear span\/} of $\bK$ and the $y_i^{(i)}$'s. 
Note that the $y_i^{(i)}$'s are linearly independent over $\bK$. 
We first define $T$ on $Y_0$. 
For $S\in \bK$ and scalars $c_1,\ldots,c_n$, set 
$$y = \sum_{j=1}^n c_j y_j^{(j)}
\tag 57$$ 
and define 
$$T(S+y) = S+P_W y\ .
\tag 58$$ 
Now if we assume that 
$$\|S+y\| =1\ ,
\tag 59$$
then 
$$|c_j| = |y_j^* (S+y)| \le M\ \text{ for all }\ j\ .
\tag 60$$
But then 
$$\|T(S+y) - (S+y)\|
= \Big\|\sum_{j=1}^n c_j P_C y_j^{(j)}\Big\| 
= 2M\ep \text{ by (56) and (59).}
\tag 61$$
Now it follows immediately that $T$ extends to a bounded linear operator 
(also denoted $T$) from $Y$ into $B(H)$, satisfying 
$$\|Ty-y\| \le 2M\ep\ \text{ for all }\ y\in Y\ .
\tag 62$$ 

Now if we assume (as we may) that $2M\ep <1$, then setting $\tilde T = TY$, 
$\tilde Y$ is a closed linear subspace of $B(H)$ and $T$ maps $Y$ one-to-one 
onto $\tilde Y$. 
Moreover, since $TY_0$ is invariant under $P_W$, so is $\tilde Y$. 
We now have that (i), (ii), (iv) of Theorem~2.1 hold (by Proposition~2.2), 
and furthermore (b) and (c) of (iii) hold. 
It remains to verify (iii)(a) and (v) of 2.1. 
Now (v) is easy, for suppose $z\in Y_0$, $z= S+y$, $S\in K$, $y$ as in (57).
Then $P_{\sim W}TZ = P_{\sim W} S+ P_{\sim W} y$, but $P_{\sim W} y$ is
actually an absolutely converging series of finite rank operators by (56). 
Thus $P_{\sim W} TZ\in K$, proving 2.1(v). 

Finally, for each $j$, define a rank-one operator $F_j:Y\to B(H)$ by 
$$F_j (y) = y_j^* (y) P_C y_j^{(j)}\ . 
\tag 63$$ 
Then it follows from (56) and (43) that 
$$\|F_j\|_{\cb} < \frac{M\ep}{2^{j-1}} \ \text{ for all }\ j\ .
\tag 64$$ 
Then setting $Q = \sum F_j$, $Q$ is of course also completely bounded, and 
$$\|Q\|_{\cb} < 2M\ep\ .
\tag 65$$
Now an inspection of the definition of $T$ on $Y_0$ yields that 
$$Tz = z-Qz\ \text{ for all }\ z\in Y\ .
\tag 66$$
But then we obtain 
$$\|T\|_{\cb} \|T^{-1}\|_{\cb} < \frac{2M\ep}{1-2M\ep}
\tag 67$$
which of course (qualitatively) yields 2.1(iii)(a).\qed
\enddemo

\remark{Remark}
Let us say that $T\in M_\infty$ is a {\it generalized block diagonal\/} 
(gbd) if there exists a $W$ of the form given in Proposition~2.2 so that 
$T = P_WT$. 
The following is a byproduct of our proof of Theorem~2.1: 
{\it Every operator in $M_\infty$ is (for every $\ep >0$) a perturbation 
of a gbd operator by a compact operator of norm less than $\ep$.} 
\endremark

\head {\smc Section} 3\endhead
\centerline{\bf The $\lambda$-Mixed Separable Extension Property}
\centerline{\bf and Extendably Locally Reflexive Banach spaces}
\bigskip

We first give the quantitative version of the property introduced in 
the preceding section. 

\proclaim{Definition} 
Let $\lambda \ge1$. 
An operator space $Z$ has the $\lambda$-Mixed Separable Extension Property 
{\rm ($\lambda$-MSEP)} if for all separable operator spaces $Y$, subspaces 
$X$, and completely bounded maps $T:X\to Z$, there exists a bounded linear 
map $\tilde T:Y\to Z$ extending $T$ with $\|\tilde T\| \le 
\lambda \|T\|_{\cb}$. 
\endproclaim

Next, we give a simple result summarizing various permanence properties 
of the MSEP for separable operator spaces $Z$. 

\proclaim{Proposition 3.1}
Let $Z$ be a separable operator space and assume $Z\subset B(H)$ for some $H$. 
Then the following are equivalent.
\roster
\item"(a)" $Z$ has the {\rm MSEP}.
\item"(b)" $Z$ is complemented in $Y$ for all separable spaces $Y$ 
with $Z\subset Y\subset B(H)$. 
\item"(c)" $Z$ has the {\rm $\lambda$-MSEP} for some $\lambda\ge1$. 
\endroster
Moreover, fixing $\lambda\ge1$, then the following are equivalent.
\roster
\item"(a$'$)" $Z$ has the {\rm $\lambda$-MSEP}.
\item"(b$'$)" $Z$ is $\lambda$-complemented in $Y$ for all $Y$ as in (b).
\item"(c$'$)" $Z$ is $\lambda$-complemented in every separable operator 
superspace.
\endroster
\endproclaim

\demo{Proof} 
(a) $\To$ (b) is essentially trivial, for let $T:Z\to Z$ be the identity map; 
a bounded linear extension $\tilde T:Y\to Z$ is a bounded projection onto $Z$.

(b) $\To$ (c) and (b$'$) $\To$ (a$'$). 
We first observe that there is a $\lambda'$ so that (b$'$) holds. 
If not, choose for every $n$, $Y_n$ a separable superspace of $Z$ contained 
in $B(H)$ so that $Z$ is not $n$-complemented in $Y_n$. 
Then letting $Y = [Y_j:j=1,2,\ldots]$, $Y$ is a separable superspace of $Z$ 
contained in $B(H)$, and $Z$ is uncomplemented in $Y$, contradicting (b). 
Now assuming  (b$'$), it suffices to show that (a$'$) holds. 

Suppose then $X,Y$ are separable operator spaces with $X\subset Y$ and 
$T:X\to Z$ is completely bounded. 
By the isometric operator injectivity of $B(H)$, we may choose $S:Y\to B(H)$  
with $\|S\|_{\cb} = \|T\|_{\cb}$ and $S|X=T$. 
Then letting $\tilde Y = [S(Y),Z]$, $\tilde Y$ is a separable superspace of $Z$
and hence there is a projection $P:\tilde Y\to Z$ onto $Z$ with $\|P\|
\le \lambda$. 
Then $\tilde T\defeq PS$ is the desired extension of $T$ with $\|\tilde T\|
\le \lambda \|T\|_{\cb}$. 
This completes the proof, in view of the triviality of the implications 
(a$'$) $\To$ (c$'$) $\To$ (b$'$).\qed
\enddemo

Although we are mainly interested in the separable case, we next note that the 
equivalence (a) $\To$ (c)  of Proposition~2.1 holds in general. 

\proclaim{Proposition 3.2} 
Let $Z$ be an operator space with the {\rm MSEP}. 
Then $Z$ has the {\rm $\lambda$-MSEP} for some  $\lambda\ge1$.
\endproclaim

\demo{Proof} 
If not, we may choose for every $n$, operator spaces $X_n$ and $Y_n$ with 
$X_n\subset Y_n$ and $T_n :X_n \to Z$ with $\|T_n\|_{\cb} =\frac1{n^2}$ 
so that 
$$\|\tilde T_n\| >n \text{ for any }  \tilde T_n :Y_n\to Z\text{ with } 
\tilde T_n|X_n = T_n\ .
\tag 68$$ 
Let now $Y = (Y_1\oplus Y_2\oplus \cdots)c_0$ and 
$X= (X_1\oplus X_2\oplus \cdots) c_0$ endowed with the standard operator 
space structure.
Of course $Y$ is separable. 
Define $T:X\to Z$ by 
$$T(x_j) = \sum T_j x_j\ .
\tag 69$$ 
$T$ is well defined, since $Z$ is a Banach space, and if $(x_j)\in X$,
then $\sum \|T_j x_j\| \le \sum \frac1{n^2} \|x\|$. 

Now given $m$, $K_1,\ldots,K_m$ in $\bK$, and $z',\ldots,z^m$ in $X$, 
we have that 
$$\align 
\Big\| \sum K_i \otimes Tz^i\Big\| 
&\le \sum_j \Big\| \sum_i K_i\otimes Tz_j^i\Big\| \tag 70\\
&\le \sum_j \|T_j\|_{\cb} \Big\| \sum_i K_i\otimes Tz_j^i\Big\|\\
&\le \sum \frac1{n^2} \max_j \Big\|\sum_i K_i\otimes Tz_j^i\Big\|\\
&= \sum \frac1{n^2} \Big\| \sum K_i z^i\Big\|\ .
\endalign$$
Hence $T$ is completely bounded, but there is no bounded linear extension 
$\tilde T:Y\to Z$ since for such a presumed extension, $\tilde T/Y_n$ 
extends $T_n$, whence $\|\tilde T|Y_n\| >n$.\qed
\enddemo

Of course the MSEP is related to a more restrictive injectivity property. 

\definition{Definition}
An operator space $Z$ is {\it mixed injective\/} if for all operator spaces 
$Y$, subspaces $X$, and completely bounded maps $T:X\to Z$, there is a 
bounded linear map $\tilde T:Y\to Z$ extending $T$. 
If, for $\lambda \ge1$, $\tilde T$ can always be chosen with 
$\|\tilde T\| \le \lambda \|T\|_{\cb}$, we say $Z$ is 
{\it $\lambda$-mixed injective\/}. 
Finally, if $Z$ is 1-mixed injective, we say that $Z$ is {\it isometrically 
mixed injective\/}. 
\enddefinition

We then have the following result, whose simple proof (via the isometric 
operator injectivity of $B(H)$) is left to the reader. 

\proclaim{Proposition} 
Let $Z$ be an operator space with $Z\subset B(H)$ for some $H$. 
Then $Z$ is mixed injective iff $Z$ is complemented in $B(H)$. 
Hence $Z$ is $\lambda$-mixed injective for some $\lambda \ge1$. 
Moreover if $\lambda\ge1$, then $Z$ is $\lambda$-mixed injective if 
$Z$ is $\lambda$-complemented in $B(H)$. 
\endproclaim 

As pointed out in the Introduction, we not not know if $\bK_0$ has the MSEP. 
The next result shows this problem is equivalent to the question of whether 
$\bK_0$ is complemented in $Y$ for all separable $Y$ with $\bK_0 \subset Y
\subset\bK_0^{**}$, in virtue of the fact that $\bK_0^{**}$ is an 
(isometrically)-injective operator space. 

\proclaim{Proposition 3.3} 
Let $Z$ be a separable operator space so that $Z^{**}$ is mixed injective. 
Then $Z$ has the {\rm MSEP} iff 
$$\text{$Z$ is complemented in $W$ for all separable spaces $W$ with }
Z\subset W\subset Z^{**}
\tag $*$ $$
\endproclaim

\demo{Proof}
One implication is trivial. 
For the slightly less trivial assertion, let $X\subset Y$ be separable 
operator spaces and $T:X\to Z$ a completely bounded map. 
Choose $\tilde T :Y\to Z^{**}$ a bounded linear extension of $\chix T$, 
where $\chix :Z\to Z^{**}$ is the canonical injection. 
Let $Y = [\chix Z,\tilde T(Y)]$ and let $P:Y\to Z$ be a surjective bounded 
linear projection (where $Z$ is of course identified with $\chix Z$). 
Then $P\tilde T$ is the desired operator extending $T$.\qed
\enddemo

\remark{Remark} 
Of course Proposition 3.2 ``reduces'' the problem of the MSEP for $\bK_0$, 
to a pure Banach space question: 
See \cite{JO} for a study of the family of separable Banach spaces $Z$ 
satisfying $(*)$, particularly in the case where $Z= (E_n)_{c_0}$, 
$E_1,E_2,\ldots$ finite-dimensional. 
\endremark

The next perhaps surprising result shows that the MSEP and mixed 
injectivity are equivalent for operator spaces complemented in their 
double duals. 

\proclaim{Proposition 3.4} 
Let $X$ be an operator space which is $\beta$-complemented in $X^{**}$ and 
suppose $X$ has the {\rm $\lambda$-MSEP}. 
Then $X$ is $\lambda\beta$-mixed injective.
\endproclaim

\proclaim{Corollary} 
Every reflexive operator space with the {\rm MSEP} is mixed injective.
\endproclaim

\demo{Proof of Proposition 3.4}
Let $Y$ be an operator super space of $X$. 
By Proposition~3.1, it suffices to prove that $X$ is 
$\lambda\beta$-complemented in $Y$.

First, fix $F$ a finite-dimensional subspace of $X$. 
We shall prove: 
$$\text{there exists a linear operator } 
T_F = T, \ T: Y\to X^{**},\ 
\text{ with } \|T\|\le\lambda
\text{ and } T|F = I|F\ .
\tag 71$$
Let $\G$ be the family of finite-dimensional subspaces of $Y$ containing 
$F$, directed by inclusion. 
For each $G\in \G$, since $X$ has the $\lambda$-MSEP, choose $T_G:G\to X$ 
a linear operator with $\|T_G\|\le\lambda$ and $T_G|F= I|F$. 
Then define $\tilde T_G:Y\to X$ by $\tilde T_G(y) = 0$ if 
$y\notin G$, $\tilde T_G(y) = T_G(y)$ if $y\in G$. 
Well, $\tilde T_G$ is neither continuous nor linear. 
However the weak*-compactness of the $\lambda$-ball of $X^{**}$ in its 
weak*-topology allows us by the Tychonoff theorem to select a subnet 
$(\tilde T_{G_\beta})_{\beta\in \D}$ 
of the net  $(\tilde T_G)_{G\in \G}$ so that 
$$\lim_{\beta\in\D} \tilde T_{\alpha_\beta} (y) \defeq T(y) 
\tag 72$$
exists weak* in $X^{**}$ for all $y\in \Ba(Y)$.
Since we do have that $\tilde T_G (\lambda y) = \lambda \tilde T_G(y)$ 
for all $y\in Y$, we obtain that the limit in (72) exists weak* for 
all $y$ in $Y$, and in fact we discover that $T$ as defined by (72) 
is indeed a linear operator with $\|T\|\le \lambda$. 
Finally, if $f\in F$, then $\tilde T_\alpha (f) =f$ for all $f$, whence 
also $Tf = f$. 
Thus (71) is proved.

Finally, let $\F$ be the family of finite-dimensional subspaces of $X$ 
directed by inclusion. 
For each $F\in \F$, choose $T_F$ satisfying (71). 
Again exploiting the weak*-compactness of the $\lambda$-ball of $X^{**}$, 
we find a subnet $(T_{F_\beta})_{\beta\in\D}$ of the net 
$(T_F)_{F\in \F}$ so that 
$$\lim_{\beta\in\D} T_{F_\beta} (y) \defeq S(y)
\tag 73$$
exists weak* in $X^{**}$ for all $y\in \Ba (y)$. 
Now it follows that $S:Y \to X^{**}$ is a linear operator 
with $\|S\|\le \lambda$. 
But if $x\in X$, then ``eventually'', $x\in F_\beta$ for $\beta\in\D$, 
whence $T_{F_\beta}(x) =x$, so also $S(x) =x$. 
Finally, letting $Q:X^{**}\to X$ be a surjective projection with $\|Q\| 
\le\beta$, we obtain that $P= QS$ is the desired projection from $Y$ 
onto $X$ of norm at most $\beta\lambda$.\qed
\enddemo

\remark{Remarks} 
1. Of course the proof shows that if $X$ is complemented in $X^{**}$, then 
$X$ is mixed injective if $X$ has the formally weaker property that 
for some $\lambda\ge1$ and for all finite-dimensional operator spaces 
$F\subset G$ and linear maps $T:F\to X$, there is a linear extension 
$\tilde T:G\to X$ with $\|\tilde T\| \le\lambda\|T\|_{\cb}$. 

2. The same compactness argument also yields that if $X$ is an operator 
space with $X$ completely complemented in $X^{**}$, then if $X$ has the 
CSEP, $X$ is injective. 
(This strengthens Proposition~2.10 of \cite{Ro2}.) 
Indeed, as noted in \cite{Ro2}, it follows that $X$ has the 
$\lambda$-CSEP for some $\lambda\ge1$. 
But then just replacing ``bounded'' by ``completely bounded'' in the 
above proof, one obtains that if $X$ is $\beta$-completely complemented 
in $X^{**}$, then $X$ is $\beta\lambda$-completely complemented in $Y$.
\endremark

We next note that certain results in \cite{R} carry over almost 
word for word to the mixed category. 

\proclaim{Proposition 3.5} 
Let $X$ be a non-reflexive operator space. 
If $X$ is mixed injective, $X$ has a subspace Banach-isomorphic to 
$\ell^\infty$. 
If $X$ is separable with the {\rm MSEP}, $X$ has a subspace Banach-isomorphic 
to $c_0$.
\endproclaim

\demo{Proof} 
If $X$ satisfies the hypothesis in the second statement, $X$ is isomorphic 
(in fact completely isomorphic) to a complemented subspace of some 
$C^*$-algebra. 
The second assertion now follows from results of H.~Pfister \cite{Pf} and 
A.~Pe{\l}czy\'nski \cite{Pe2}. 
If $X$ satisfies the first hypothesis, $X$ is completely isomorphic to a 
complemented subspace of some von-Neumann algebra. 
The first assertion now follows from these results and the result 
of \cite{Ro1} that any non-weakly compact operator from $\ell^\infty$ into 
some Banach space fixes a copy of $\ell^\infty$. 
The argument itself is word for word as the proof of Proposition~2.8 of 
\cite{Ro2}, deleting the word ``completely'' in all its 
occurrences.\qed
\enddemo

Finally, we note the analogue of Proposition 2.22 of \cite{Ro2}. 

\proclaim{Proposition 3.6} 
Let $X$ be a separable operator space with the {\rm $\lambda$-MSEP}. 
If $\lambda <2$, then $X$ is reflexive (and hence is $\lambda$-mixed 
injective by Proposition~3.3). 
\endproclaim

\demo{Proof} 
The argument is essentially the same as that for Proposition 2.22 of 
\cite{Ro2}. 
We give this argument however, for the sake of completeness. 
Suppose to the contrary that $X$ is not reflexive. 
Then $X$ contains a subspace isomorphic to $c_0$ by Proposition~3.5. 
Now let $\ep>0$, to be decided later, and choose (using the ``folklore'' 
result, proved in Proposition~2.22 of \cite{Ro2}) a subspace $Z$ of $X$ 
which is Banach $(1+\ep)$-isomorphic to $c_0$ and $(1+\ep)$-complemented 
in $X$. 
Now let $Y$ be a separable subspace of $Z^{**}$ with $Z\subset Y$ and let  
$i:Z\to X$ be the identity injection, and also let $P:X\to Z$ be a 
surjective projection with $\|P\| <1+\ep$. 
Since  $X$ has the $\lambda$-MSEP, letting $Y$ have its natural operator 
space structure, we find a bounded linear extension $\tilde\imath:Y\to X$ 
with $\|\tilde\imath\| \le\lambda$. 
But then letting $Q= P\tilde\imath$, $Q$ is a projection from $Y$ 
onto $Z$ and 
$$\|Q\| < (1+\ep) \lambda\ .
\tag 74$$
Since $Z$ is $(1+\ep)$-isomorphic to $c_0$, it now follows that if 
$\tilde Y$ is separable with $c_0\subset\tilde Y\subset \ell^\infty$, then 
$$c_0\text{ is $(1+\ep)^2\lambda$-complemented in } \tilde Y\ .
\tag 75$$
This implies $c_0$ itself has the $(1+\ep)^2\lambda$ SEP, whence by a 
result of Sobczyk \cite{S}, $(1+\ep)^2\lambda\ge2$. 
Of course this is a contradiction for $\ep$ small enough.\qed
\enddemo

We now give some examples of operator spaces with the MSEP. 
Evidently any complemented subspace of an operator space with the MSEP 
also has the MSEP. 
The next result is thus an immediate consequence of a result in \cite{Ro2}. 

\proclaim{Proposition 3.7} 
Let $X$ be a $\lambda$-complemented subspace of $c_0 (R\oplus C)$. 
Then $X$ has the {\rm $2\lambda$-MSEP}. 
\endproclaim

\demo{Proof} 
It is proved in \cite{Ro2} that $c_0(R\oplus C)$ has the 
2-CSEP, hence trivially the 2-MSEP.\qed
\enddemo

Of course $c_0(R\oplus C)$ is Banach isomorphic to $c_0(\ell^2)$, and the 
infinite-dimensional complemented subspaces of $c_0(\ell^2)$ have 
been isomorphically classified in \cite{BCLT}; there are exactly six of them.

\proclaim{Problem 3.1} 
Let $X$ be an infinite-dimensional separable operator space with the 
{\rm MSEP}. 
Is $X$ Banach isomorphic to one of the spaces 
$$c_0,\quad (\ell_n^2)_{c_0},\quad c_0(\ell^2),\quad \ell^2,\quad 
c_0\oplus \ell^2\ ,\ \text{ or }\ (\ell_n^2)_{c_0}\oplus \ell^2\ ?
\tag 76$$
\endproclaim

Of course if $\bK$ (or equivalently $\bK_0$) {\it has\/} the MSEP, the 
answer is negative, and the list must be much bigger than this.
It is worth pointing out, however, that work of G.~Pisier yields 
immediately that 
{\it every reflexive mixed injective operator space is Hilbertian\/}, i.e., 
Banach isomorphic to a Hilbert space (cf.\ \cite{R}). 
Thus the list in (76) {\it is\/} complete in the separable reflexive case; 
$\ell^2$ is the only example!

We next give some examples of 1-mixed injective operator spaces. 
Let $\nat^*= \nat\cup  \{\infty\}$ and $n,m\in\nat^*$. 
Recall that $\M_\infty$ denotes $B(\ell^2)$ regarded as matrices 
operating on the natural basis. 

Now the following are all 1-mixed injective.
\roster
\item"I." $\M_{n,m}$
\item"II." $ \SS_n$, the $n\times n$ symmetric matrices
\item"III." $\A\SS_n$, the anti-symmetric $n\times n$ matrices.
\endroster
(If $A^t$ denotes the transpose of $A$, then $A\in \SS_n$ iff 
$A= A^t$; $A\in \A\SS_n$ iff $A= -A^t$.) 

Neither $\SS_\infty$ nor $\A\SS_\infty$ are injective, while of course 
$\M_{n,m}$ is 1-injective for all $n$ and  $m$. 
However {\it another\/} family of 1-mixed injectives occurs; the 
{\it spin factors\/}. 

\definition{Definition}
A closed subspace $X$ of $B(H)$ is called a {\it spin factor\/} if 
\roster
\item"(a)" $X$ {\it is self-adjoint}
\item"(b)" $\dim X>1$
\item"(c)" the square of every element of $X$ is a scalar.
\endroster
\enddefinition

It is known that spin-factors are 1-mixed injective \cite{ES} 
and Hilbertian. 
Moreover, in the separable case, $X$ is a spin-factor iff there exists a 
sequence $S_1,S_2,\ldots$ of anti-commuting self-adjoint unitaries with 
$X= [S_n]$. 
Here, $(S_n)$ is either finite of length at least 2, or infinite. 
$X$, as an operator space, is determined up to complete isometry by its 
dimension (i.e., the length of this sequence $(S_n)$). 
For $n\in\nat^* \sim \{1\}$, let $\Sp (n)$ denote a spin factor of dimension 
$n$ for $n<\infty$ (resp. separable infinite-dimensional if $n=\infty$).

Standard constructions yield that for all $n$, $\Sp (n)$ is 1-completely 
isometric to a (necessarily contractively complemented) subspace of 
$\M_{2^{n/2}}$ if $n$ is even, $\M_{2^{[n/2]}} \oplus \M_{2^{[n/2]}}$ if 
$n$ is odd. 

However the following result yields that $\Sp (\infty)$ is not completely 
isomorphic to a subspace of $K$. 

\proclaim{Proposition 3.8} 
Let $X$ be an operator space so that $X\otimes_{\op}X$ is completely 
isomorphic to a subspace of $X$. 
If $\Sp (\infty)$ completely embeds in $X$, then $\ell^1$ Banach embeds 
in $X$, hence $X^*$ is non-separable. 
\endproclaim

\demo{Proof} 
By a result of U.~Haagerup \cite{H} (see also \cite{Pa}), 
if $(S_i)$ is an infinite spin system 
in $B(H)$, $(S_i\otimes S_i)$ is Banach-equivalent to the usual 
$\ell^1$-basis. 
Now if $Y$ is a subspace of $X$ which is completely isomorphic to 
$\Sp(\infty)$, $Y\otimes_{\op} Y$ is completely isomorphic to 
$\Sp (\infty) \otimes_{\op} \Sp (\infty)$, hence $\ell^1$ embeds 
in $X\otimes_{\op}X$ by Paulsen's results.\qed
\enddemo 

We next give a remarkable, simple permanence property of mixed injectivity 
and the Mixed Separable Extension Property. 
We first need the following (apparently new) concept. 

\definition{Definition}
Given operator spaces $X$ and $Y$, $X$ is {\it completely semi-isomorphic\/} 
to $Y$ 
if there exists a completely bounded map $T:X\to Y$ which is a Banach 
isomorphism, i.e., (since $X,Y$ are assumed complete), so that $T$ is 1--1
and onto. 
We call such a map $T$ a {\it complete semi-isomorphism\/} from $X$ onto $Y$. 
In case $\|T\|_{\cb} =1 = \|T^{-1}\|$, $T$ is called a {\it complete 
semi-isometry\/} and $X$ is said to be {\it completely semi-isometric\/} to 
$Y$ in case there exists a complete semi-isometry mapping $X$ onto $Y$. 
In case $\|T\|_{\cb} \|T^{-1}\| \le\lambda$, we say $X$ is 
$\lambda$-completely semi-isomorphic to $Y$.  
Finally, we set $d_s (X,Y) = \inf \{\lambda :X$ is completely semi-isomorphic 
to $Y\}$. 
\enddefinition

This relation is trivially reflexive and is also (quantitatively) 
transitive: 
if $X$ is $\lambda$-completely semi-isomorphic to $Y$ and $Y$ is 
$\beta$-completely semi-isomorphic to $Z$, 
then $X$ is $\lambda\beta$-completely semi-isomorphic to $Z$. 
The relation is of course not symmetric in general; e.g.,  
$R\cap C$ is completely semi-isomorphic to $R$ but $R$ is not 
completely semi-isomorphic to $R\cap C$. 
It can be shown that the relation does not yield 
a partial order on operator spaces modulo complete isomorphism. 
In fact, there exist non-completely isomorphic operator spaces $X$ and 
$Y$ so that each is completely semi-isometric to the other. 
However if $X$ and $Y$ are each completely semi-isomorphic to the other, 
then $X$ and $Y$ {\it are\/} completely isomorphic 
if one of them, say $X$, is homogeneous, i.e., if every bounded 
operator on $X$ is completely bounded. 
Indeed, suppose $X$ is $\lambda$-homogeneous 
(i.e., $\|U\|_{\cb} \le \lambda \|U\|$ for all $U\in \L(X)$) 
and suppose $T:X\to Y$ 
and $S:Y\to X$ are surjective complete semi-isomorphisms. 
But then $T^{-1} S^{-1}$ is completely bounded, hence so is $T^{-1} = 
T^{-1} S^{-1} S$, and $\|T^{-1}\|_{\cb} \le \lambda \| T^{-1}\|\, \|S^{-1}\| 
\, \|S\|_{\cb}$. 
We thus obtain that $d_{\cb}(X,Y) \le \lambda d_s (X,Y)d_s (Y,X)$. 

We now give the permanence property mentioned above: 
mixed injectivity and the MSEP are both preserved under 
complete semi-isomorphisms; 
i.e., if $X$ is completely semi-isomorphic 
to $Y$ and $Y$ has one of these properties, so does $X$. 

\proclaim{Proposition 3.9} 
Let $\lambda,\beta\ge1$ and let $Z$ and $\tilde Z$ be operator spaces with 
$\tilde Z$ $\beta$-completely semi-isomorphic to $Z$. 
Then if $Z$ is $\lambda$-mixed injective (resp. has the 
{\rm $\lambda$-MSEP\/}), $\tilde Z$ is $\beta\lambda$-mixed injective 
(resp. has the {\rm $\lambda\beta$-MSEP\/}). 
\endproclaim

\demo{Proof}
Choose $S:\tilde Z\to Z$ a surjective complete semi-isomorphism with 
$\|S\|_{\cb} \|S^{-1}\|\le \beta$. 
Let $X\subset Y$ be operator spaces and $T:X\to \tilde Z$ be a completely 
bounded map. 
Now consider the diagram:
$$\lower.5truein\hbox{\epsfysize=1.0truein\epsfbox{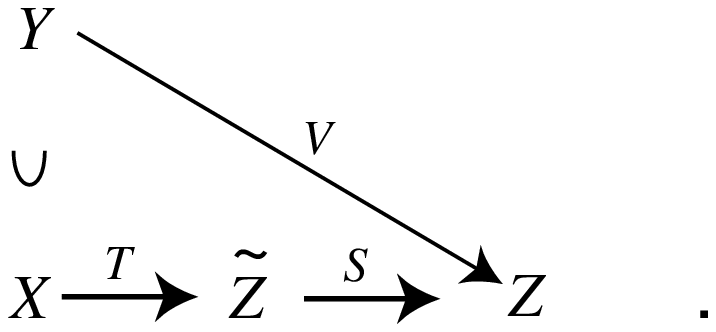}}$$
In the case where $Z$ is $\lambda$-mixed injective, choose $V$ a linear 
operator completing this diagram  with 
$$\|V\| \le \|S\circ T\|_{\cb} \le \lambda \|T\|_{\cb} \|S\|_{\cb}\ .
\tag 77$$
Then letting $\tilde T= S^{-1}V$, we obtain that $\tilde T:Y\to\tilde Z$ 
is an extension of $T$ with 
$$\align
\|\tilde T\| & \le \lambda\|S\|_{\cb} \|S^{-1}\|\, \|T\|_{\cb}
\quad\text{ by (77)}\tag 78\\
&\le \lambda\beta\|T\|_{\cb}\ .
\endalign$$
Of course (78) yields $\tilde Z$ is $\lambda\beta$-mixed injective. 

In the case where $Z$ has the $\lambda$-MSEP, simply assume that $Y$ 
is in addition separable, to obtain the desired conclusion.\qed
\enddemo

Proposition 3.9 has the immediate consequence: 
If $X$ is completely semi-isomorphic to a space with the CSEP, 
$X$ has the  MSEP. 
This suggests the following question. 

\proclaim{Problem 3.2} 
Let $X$ be a separable operator space with the {\rm MSEP}. 
Is $X$ completely semi-isomorphic to a space with the {\rm CSEP}? 
\endproclaim

We do not know if $\bK$ is semi-isomorphic to a space with the CSEP, 
although we suspect this is not the case. 
Let us note, however, that the presently known examples of separable 
spaces with the MSEP do have the property specified in this problem, 
e.g., $\Sp (\infty)$ is completely semi-isomorphic to $R$.

The 1-mixed injective  finite-dimensional spaces are completely classified 
up to Banach isometry, 
based in part on deep work of E.~Cartan \cite{C} for which there seems 
to be no decent modern exposition. 

\proclaim{Theorem A} 
Let $X$ be a finite-dimensional isometrically injective operator space. 
Then $X$ is Banach isometric to a (finite) $\ell^\infty$-direct sum of 
spaces $E$ each of the following form for some $m,n\in \nat$
\roster
\item"I." $E= \M_{n,m}$
\item"II." $E= \SS_n$ 
\item"III." $E= \A\SS_n$
\item"IV." $E= \Sp(n)$
\endroster
\endproclaim
Spaces of the form I--IV are known as {\it Cartan\/} factors of types I--IV.
Of course any $\ell^\infty$ finite direct sum of Cartan factors of 
types I--IV {\it is\/} isometrically mixed injective (in its natural 
operator space structure). 
Now Proposition~3.8 coupled with the spaces listed in Theorem~A, yields 
a rather vast supply of finite-dimensional 1-mixed injectives 
(e.g., $R_n\cap C_n$ is of this form). 
Are these the only ones?

\proclaim{Problem 3.3} 
Is every finite dimensional 1-mixed injective completely semi-isometric to an 
$\ell^\infty$ direct sum of Cartan factors of types I--IV?
\endproclaim

The work in \cite{AF} is certainly related to this problem, especially 
if the answer is negative!
Problem~3.2 ought to be solved in this century! 
The next problem, on the other hand, seems quite intractable at this time
(although a negative answer need not be). 
An affirmative solution would imply an affirmative solution to the famous 
``finite-dimensional ${\Cal P}_\lambda$ problem'' in the commutative theory.

\proclaim{Problem 3.4} 
Given $\lambda>1$, is there a $\beta$ so that every $\lambda$-mixed injective 
finite-dimensional space is $\beta$-completely 
semi-isomorphic to a 1-mixed injective?
\endproclaim

A remarkable factorization theorem due to M.~Junge \cite{J} yields that a 
purely local formulation of the classification problem for finite-dimensional 
mixed injectives. 
We are indebted to M.~Junge for the proof of this result, which yields that 
the finite-dimensional $\beta$-mixed injectives are essentially, up to 
complete semi-isomorphism, the $\beta$-complemented subspaces of $\M_n$'s. 

\proclaim{Proposition  3.10} 
Let $X$ be a finite-dimensional operator space and $\lambda\ge1$. 
The following are equivalent:
\roster
\item"(1)" $X$ is $\lambda$-mixed injective.
\item"(2)" For all $\ep>0$, there exist an $n$ and linear maps 
$U:X\to \M_n$ and $V:\M_n\to X$ so that $I_X = VU$ and $\|V\|\, \|U\|_{\cb} 
<\lambda +\ep$. 
That is, we have the diagram 
$$\matrix 
\M_n\\
\llap {\raise.5ex\hbox{$\scriptstyle U$}}\!\!\!\nearrow\qquad
\searrow \rlap{\raise.5ex\hbox{$\!\!\!\scriptstyle V$}}\\
X\quad  \buildrel I\over \longrightarrow\quad  X\ .
\endmatrix$$
\endroster
\endproclaim

\proclaim{Corollary 3.11} 
If $X$ is finite-dimensional and $\lambda$-mixed injective, then for all 
$\ep>0$, there is a subspace $Y$ of $\M_n$ so that $X$ is  
$(\lambda+\ep)$-completely 
semi-isomorphic to $Y$ and $Y$ is $(\lambda+\ep)$-mixed 
injective. 
\endproclaim
\endremark

\remark{Remark} 
Of course the conclusion of the Corollary implies that $Y$ is 
$(\lambda +\ep)$-Banach complemented in $\M_n$.
\endremark

\demo{Proof of 3.11} 
Set $Y= U(X)$, where $U,V$ are chosen as in (2) of 3.10. 
It follows that $U:X\to Y$ is a semi-isomorphism with $U^{-1}= V|_Y$, 
hence $d_s (X,Y) \le \|V|_Y\|\, \|U\|_{\cb} < \lambda +\ep$. 
Setting $P= UV$, then $P$ is a projection from $\M_n$ onto $Y$, 
and $\|P\| \le \|U\|\, \|V\| <\lambda +\ep$, as desired.\qed
\enddemo 

\demo{Proof of Proposition 3.10}
(1) $\To$ (2): 
Assume without loss of generality that $X\subset B(H)$. 
Then there exists a surjective linear projection $P:B(H)\to X$ with 
$\|P\| \le\lambda$. 
Let $Y= (X,\MIN)$, and let $T:X\to Y$ be the formal identity map  and 
$i:X\to B(H)$ be the identity injection. 
Thus $T$ completely factors  through $B(H)$, $T =Pi$, and 
$\|P\|_{\cb} \|i\|_{\cb} \le\lambda$. 
Hence by a basic factorization theorem in \cite{J} 
(reproved as Theorem~7.6 in \cite{EJR}; see also Remark~3.6 in \cite{JM}), 
we may choose $n$ and linear maps 
$U:X\to \M_n$, $\tilde V:\M_n \to Y$ with $T= \tilde VU$ and 
$\|\tilde V\|_{\cb}\|U\|_{\cb} < \|P\|_{\cb} \|i\|_{\cb} +\ep 
\le \lambda+\ep$. 
But now just letting $V= T^{-1} \tilde V$, we obtain (2) of 3.10, for 
trivially $\|V\| \le \|T^{-1}\|\, \|\tilde V\|_{\cb} = \|\tilde V\|_{\cb}$. 

(2) $\To$ (1):
Let $Y\subset Z$ be separable operator spaces and $T:Y\to X$ be a given 
linear map. 
Let $\ep >0$ and choose $\M_n$ and $U,V$ as in (2). 
Now since $\M_n$ is 1-injective, choose $S:Z\to \M_n$ a linear map with 
$\|S\|_{\cb} = \|UT\|_{\cb} \le \|U\|_{\cb}\|T\|_{\cb}$. 
Thus we have the diagram
$$\gather \qquad\quad \M_n\\
\hbox{\epsfysize=0.65truein\epsfbox{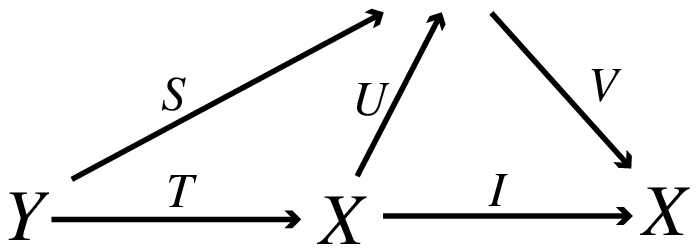}}
\endgather$$
Then $\tilde T_\ep  \defeq VS$ extends $T$ and $\|\tilde T_\ep\|_{\cb} 
< \lambda +\ep$. 
Since $X$ is finite-dimensional, we may choose a sequence $(\ep_n)$ 
tending to zero and an operator $\tilde T:Z\to X$ so that 
$\tilde T_{\ep_n} \to \tilde T$ in the strong operator topology. 
It follows that $\tilde T$ extends $T$ and $\|\tilde T\|_{\cb} \le\lambda$. 
Thus $X$ has the $\lambda$-MSEP, so by Proposition~3.4, $X$ is 
$\lambda$-mixed injective.\qed
\enddemo

We briefly indicate the remarkable connection of Theorem~A with a rather vast 
domain of modern research. 
A closed linear subspace $X$ of $B(H)$ is called a (concrete) 
$JC^*$-triple if 
$$TT^*T \in X\ \text{ whenever }\ T\in X\ .
\tag 79$$
It then follows by polarization that 
$$2\{A,B,C\} \defeq AB^* C + CB^*A\text{ belongs to $X$ whenever $A,B,C$ do}.
\tag 80$$
$\{A,B,C\}$ is called the {\it triple product\/} of $A$, $B$, and $C$; 
an abstract generalization of this led to the theory of $JB^*$ triples.
In turn, this theory yields the following remarkable general result, which 
includes Theorem~A. 

\proclaim{Theorem B} 
Let $X$ be a finite-dimensional Banach space.
The following are equivalent.
\roster
\item"1." $X$ is isometric to a contractively complemented subspace of some
$C^*$-algebra.
\item"2." $X$ is isometric to a $JC^*$-triple.
\item"3." The open unit ball of $X$ is biholomorphically transitive, and $X$ 
contains no contractively complemented subspace isometric to a Cartan factor 
of type~V.
\item"4." $X$ is of the form specified in Theorem~A.
\endroster
\endproclaim

Several of the implications in Theorem B hold in infinite-dimensions as well. 
In fact, it is known that $1\To 2\To 3$ in general and finally if $X$ satisfies 
3.\ and is isometric to a dual space, then $3\To1$. 
$1\To 2$ is due to Y.~Friedman and B.~Russo \cite{FR} and 
$2\To3$ (without the assertion concerning the type~V factor, which came later) 
is due to L.~Harris \cite{H}; see also \cite{K2}. 
For $3\To1$ for $X$ a dual, see C.-H.~Chu and B.~Iochum \cite{CI}. 
As far as we know, the following are open questions in general: 
Does $3\To1$? 
Does $2\To1$?
The profound result of Cartan's which underlies this:  
{\it the unit ball of a
finite-dimensional Banach space is biholomorphically transitive iff the 
space is an $\ell^\infty$-direct sum of Cartan factors\/}. 
In addition to the factors of types~I--IV, there are two more, 
types~V and VI; the type~VI factor 
consists of the $3\times 3$ Hermitian matrices 
over the  complex octonions, and is 27-dimensional. 
The type~V factor embeds in this one; as a Banach space, it may however be 
explicitly identified as follows:
Let $e_0,e_1,\ldots,e_7$ be the usual basis for the complex octonions $O$ 
(with $e_0$ the identity).  
For $a\in O$, $a= \sum_{i=0}^7 a_i e_i$ with the $a_i$'s complex scalars, 
set $|a| = (\sum |a_i|^2)^{1/2}$ and $n(a) = \sum a_i^2$; 
also set $\tilde a = a_0 e_0 - \sum_{i=1}^7 a_i e_i$. 
Now the Cartan factor of type~V may be identified with $X= O\times O$ where, 
if $x = (a,b)$, then 
$$\|x\|^2 = |a|^2 + |b|^2 + 
\sqrt{(|a|^2 + |b|^2)^2 - |n(a)|^2 + |n(b)|^2 + 2|\tilde a\,b|^2}\ .
\tag 81$$
Note that if $a$ and $b$ have only real coefficients, $|\tilde a\,b| = 
|\tilde a|\,|b| = |a|\,|b|$ by a fundamental property of the real octonions, 
whence 
\roster
\item"{}" 
$\|x\|^2 = |a|^2 + |b|^2$, the ordinary Euclidean norm of the vector $x$.
\endroster
(The  industrious reader may dig the proof of Theorem~B out of the 
references \cite{C}, \cite{FR}, \cite{K1}, \cite{K2}. 
See also \cite{H} for important earlier structure results on $JC^*$-triples. 
Also see \cite{CI} and \cite{LO}, \cite{D}, and finally \cite{Ru} for 
a comprehensive survey on $JB^*$-triples. 
Also, although (81) is a simple deduction from known work, this explicit 
expression for the actual {\it norm\/} on the type~V Cartan factor, seems 
to be new.) 

To further penetrate the fundamental question of whether $\bK$ has the 
MSEP, we introduce the following new concept in pure Banach space theory.

\definition{Definition}
A Banach space $X$ is called {\it Extendably Locally Reflexive\/} (ELR) 
if there exists a $\lambda\ge1$ so that for all finite dimensional subspaces 
$F$ and $G$ of $X^*$ and $X^{**}$ respectively and all $\ep>0$, there exists 
an operator $T:X^{**}\to X^{**}$ with 
$$\left\{ \eqalign{
\text{(i)}\hfil\quad&TG\subset X\cr
\text{(ii)}\hfil\quad&\langle Tg,f\rangle = \langle g,f\rangle \text{ for all } 
g\in G,\ f\in F\cr
\text{(iii)}\hfil\quad&\|T\| < \lambda +\ep\ .\cr}\right.
\tag 82$$
In case  $\lambda$ works, we say $X$ is $\lambda$-ELR. 
\enddefinition

The terminology is motivated as follows: 
by the Local Reflexivity Principle due jointly to J.~Lindenstrauss and the 
second author of this paper \cite{LR2} (see also \cite{JRZ}); 
{\it for all $X$, and $F,G$ as above, $\ep>0$, 
there exists an operator $T:G\to X$ with $\|T\| <1+\ep$ and 
satisfying\/} (82)(ii). 
Then  $X$ is ELR precisely when 
there exist such operators which admit uniformly 
bounded extensions $\tilde T$ to all of $X^{**}$, i.e., we have 
$$\CD
X^{**}@>{\tilde T}>> X^{**}\\
\bigcup @. \bigcup\\
G @>T>> X
\endCD 
\quad\text{with }\ \|\tilde T\|\le C$$
for some absolute constant $C$. 

The next result yields several equivalences for Extendable Local 
Reflexivity. 

\proclaim{Proposition 3.12} 
Let $\lambda \ge1$, $X$ a given Banach space. 
The following are equivalent:
\roster
\item"(i)" $X$ is $\lambda$-{\rm ELR}.
\item"(ii)" there exists a net $(T_\alpha)$ of linear operators on $X^{**}$ 
with $\|T_\alpha\|\le\lambda$ for all $\alpha$, so that for all $x^{**}\in 
X^{**}$
\item"{}" \quad {\rm (a)} $T_\alpha x^{**}\to x^{**}$ weak*
\item"{}" and
\item"{}" \quad {\rm (b)} $T_\alpha x^{**}$ is ultimately in $X$.
\item"(iii)" same as (ii), with the addition 
\item"{}" \quad {\rm (c)} $T_\alpha x\to x$ in norm, for all $x\in X$.
\item"(iv)" for all $F,G$ finite-dimensional subspaces of $X^*$ and $X^{**}$ 
respectively, there is an operator $T:X^{**}\to X^{**}$ satisfying 
{\rm (82)(i)--(iii)} and in addition
$$Tg =g\text{ for all } g\in G\cap X\ .
\tag 83$$
\endroster
\endproclaim

\demo{Proof} 
(i) $\To$ (ii). 
Let $\D = \{F,G,\ep :F,G$ are finite-dimensional subspaces of $X^*$ and 
$X^{**}$ respectively, and $0<\ep <1$. 
Direct $\D$ by:
\roster
\item"{}" $(F,G,\ep) \le (F',G',\ep')$ if $F\subset F'$, $G\subset G'$, and 
$\ep' \le \ep$.
\endroster
Given $\alpha = (F,G,\ep)$ in $\D$, choose $T\alpha =T$ satisfying (82), 
and set $T_\alpha = \frac{\lambda}{\lambda+\ep} T_\alpha$. 
Then $(T_\alpha)_{\alpha\in\D}$ has the desired property. 

(ii) $\To$ (iii). 
Let $\D$ be the directed set given in the above proof, and also suppose 
the net satisfying (ii) is given by $(T_\alpha)_{\alpha\in\G}$.
Now given $d= (F,G,\ep)$ in $\D$, choose $\beta\in\G$ so that for all 
$\alpha\ge\beta$, 
$$\eqalign{
\text{(i)}\hfil\quad&T_\alpha G\subset X\cr
\text{(ii)}\hfil\quad&|\langle T_\alpha g,f\rangle - \langle g,f\rangle|
\le \ep\| g\|\, \|f\|\text{ for all } g\in G,\ f\in F\ .\cr}
\tag 84$$
Now if $x\in G\cap X$, then $T_\alpha x\to x$ weakly; hence certain 
far out convex combinations converge in norm. 
But then, thanks to the finite-dimensionality of $G$, we may choose a convex 
combination $S_d$ of $\{T_\alpha: \alpha\ge\beta\}$ so that 
$$\|S_d x-x\| \le \ep \|x\|\text{ for all } x\in G\cap X\ .
\tag 85$$

Now $S_d$ still satisfies (84)(i) (replacing ``$T_\alpha$'' by ``$S_d$'' there)
and of course $\|S_d\| \le\lambda$ also, 
hence it follows that the net $(S_d)_{d\in\D}$, satisfies the 
conclusion of (iii). 

(iv) $\To$ (i) -- trivial. 

(iii) $\To$ (iv). 
Let $\ep>0$ and fix $F,G$ finite-dimensional subspaces of $X^*$, $X^{**}$ 
respectively, choose $f_1,\ldots,f_n$ a basis for $F$, and choose 
$x_1,\ldots,x_n$ in $X$ with $f_i (x_j) = \delta_{ij}$ for all $i$ and $j$. 
Now assuming $(T_\alpha)$ satisfies (iii), we may choose $\alpha$ so that 
(84)(i) holds and also 
$$\eqalign{
\text{(i)}\hfil\quad&|\langle T_\alpha g-g),f_i\rangle|\le \delta\|g\| 
\text{ for all } g\in G\cr
\text{(ii)}\hfil\quad&\|T_\alpha g-g\| \le \ep \|g\|\text{ for all } 
g\in G\cap X\ .\cr}
\tag 86$$

Now also choose $Y$ a linear subspace of $G$ with 
$$Y\oplus (G\cap X) = G\ ;
\tag 87$$
then choose $P,Q$ linear projections onto $Y$, $G\cap X$ respectively so that 
$$G\cap X \subset \ker P\text{ and } Y\subset \ker Q\ .
\tag 88$$
Of course then by (87) and (88), $QP = PQ=0$, and 
$$P|G\text{ is a projection onto $Y$ with kernel } G\cap X,\ \text{ and } 
Q|G= (I-P)|G\ .
\tag 89$$
Let $\delta>0$, to be decided later. 
Assuming $(T_\alpha)$ satisfies (iii) of the Theorem, choose $\alpha$ so 
that (84)(i) holds and also 
$$\eqalign{
\text{(i)}\hfil\quad&|\langle (T_\alpha g-g),f_i\rangle| \le \delta \|g\|
\text{ for all } g\in G\cr 
\text{(ii)}\hfil\quad&\|T_\alpha g-g\| \le \delta \|g\|\text{ for all } 
g\in G\cap X\ .\cr}
\tag 90$$
Then define $T:X^{**}\to X^{**}$  by 
$$Tz = \sum_{j=1}^n \langle Pz-T_\alpha Pz,f_j\rangle x_j + 
T_\alpha Pz + Qz + T_\alpha Rz\
\tag 91$$
for all $z\in X^{**}$.
Now if $z\in G\cap X$, then $P(z) = R(z) =0$ and $Q(z) =z$, so $Tz=z$; 
hence (83) holds, and so $$\langle Tz,f\rangle = \langle z,f\rangle
\text{ for all } f\in F\ .
\tag 92$$
If $z\in Y$, then $Q(z) = R(z)=0$ and $P(z) = z$; whence 
$$Tz = \sum_{j=1}^n \langle z,f_j\rangle x_j 
- \sum_{j=1}^n \langle T_\alpha  z,f_j\rangle x_j +T_\alpha z\ .
\tag 93$$
But then for each $j$, 
$$\langle Tz,f_j\rangle = \langle z,f_j\rangle - \langle T_\alpha x,f_j
\rangle + \langle T_\alpha z,f_j\rangle = \langle z,f_j\rangle\ .
\tag 94$$
Since the $f_j$'s are a basis for $F$, (92) holds. 
But then since (93) holds for $z\in G\cap X$ and $z\in Y$, (82)(ii) holds. 
Finally, we estimate the norm of $T$. 
Now fixing $z$ in $X^{**}$, $z=Pz+Qz+Rz$. 
Hence 
$$(T-T_\alpha)z = \sum_{j=1}^n \langle Pz-T_\alpha Pz,f_j\rangle x_j 
+ Qz-T_\alpha Qz\ .
\tag 95$$
Thus, we obtain by (90) that 
$$\|(T-T_\alpha) z\| = \delta \|P\| \sum_{j=1}^n \|x_j\|\, \|z\| + 
\delta \| Q\|\, \|z\|\ .
\tag 96$$
Hence, simply choosing $\delta$ so small that 
$$\delta \biggl( \|P\| \sum_{j=1}^n \|x_j\| + \|Q\|\biggr) <\ep\ ,
\tag 97$$
we obtain that the finite-rank perturbation $T-T_\alpha$ of $T$ has norm 
smaller than $\alpha$, whence 
$$\|T\| < \|T_\alpha \| + \ep \le \lambda +\ep\ . 
\tag 98$$
completing the proof of Proposition 3.12.\qed
\enddemo

Although not evident from the definition  of the Extendable Local 
Reflexivity, there is an astonishing connection between this property 
and the bounded approximation property (the bap).

\proclaim{Theorem 3.13}
Let $X$ be a given Banach space. 
The following assertions  are equivalent.
\roster
\item"(i)" $X$ is {\rm ELR} and has the {\rm bap}.
\item"(ii)" There exists a uniformly bounded net $(T_\alpha)$ of finite 
rank operators on $X^{**}$ with $T_\alpha x^{**}\to x^{**}$ weak* 
for all $x^{**}\in X^{**}$.
\item"(iii)" $X^*$ has the {\rm bap}.
\endroster
\endproclaim

\remark{Remarks}

1. The second author of this paper discovered the ELR concept as well as 
the implication (iii) $\To$ (i)  during a research visit to Odense 
University, November 1997. 
(Of course $X^*$ has bap $\To X$ has bap is an old standard result.) 
The implication (i) $\To$ (iii) was discovered by the authors of \cite{JO} 
shortly after an initial draft of the present paper was prepared. 

2. Our proof yields that one can choose a net $(T_\alpha)$ satisfying (ii) 
with $\|T_\alpha\|\le \lambda$ for all $\alpha$ iff $X^*$ has the 
$\lambda$-bap. 
On the other hand, if $X$ is $\lambda$-ELR and has the $\beta$-bap, $X^*$ has 
the $\lambda\beta$-bap (as also obtained in \cite{JO}). 

3. It follows immediately from Theorem~3.13 
and a deep result of T.~Szankowski \cite{S} that $C_1$ {\it fails\/} to be 
ELR ($C_1$ the trace class operators on Hilbert space). 
This and other examples of Banach spaces failing to be ELR are given 
in \cite{JO}.
\endremark

\demo{Proof of Theorem 3.13}
(i) $\To$ (ii). 
Suppose $X$ is $\lambda$-ELR and has the $\beta$-bap. 
Let $\D$ be the directed set given in the proof of (i) $\To$ (ii) of the 
preceding Proposition.
Given $\alpha = (F,G,\ep)$ in $\D$, first choose $T:X^{**}\to X^{**}$ 
satisfying (82). 
Now choose $S:X\to X$ a finite rank operator with 
$$\|S\| < \beta +\ep \text{ and } Sx=x\text{ for all } x\in TG\ .
\tag 99$$
Hence 
$$\langle S^{**}Tg,f\rangle  = \langle STg,f\rangle 
= \langle g,f\rangle\text{ for all } g\in G,\ f\in F,
\text{ by (82) and (99)}\ .
\tag 100$$
Finally, let $U_\alpha = S^{**}T$. 
Then $U_\alpha x^{**}\to x^{**}$ $\omega^*$ for all $x^{**}\in X^{**}$, 
$$\|U_\alpha\|\le (\lambda +\ep)(\beta+\ep)$$
for all $\alpha$, 
and moreover $\varlimsup_{\alpha\in\D} \|U_\alpha\| \le\lambda\beta$. 
So if we let $T_\alpha = \frac{\lambda\beta}{(\alpha+\ep )(\beta+\ep)}
U_\alpha$, then $(T_\alpha)$ satisfies (ii) with $\|T_\alpha\| \le 
\lambda\beta$ for all $\alpha$ (as claimed in Remark~1 above). 

(ii) $\To$ (iii). 
Let $(T_\alpha)$ be a net of finite rank operators satisfying (ii), and 
suppose $\|T_\alpha\| \le\lambda$ for all $\alpha$. 
Now since the $T_\alpha$'s are finite rank, it follows that we may assume 
the $T_\alpha$'s are {\it weak*-continuous\/}. 
To see this, again let $\D$ be the directed set given above. 
For $\beta\in\D$, $\beta = (F,G,\ep)$, choose $\alpha$ so that 
$$|\langle T_\alpha g,f\rangle - \langle g,f\rangle| 
\le (1+\ep)\|g\|\, \|f\| 
\tag 101$$
for all $g\in G$ and $f\in F$. 

\noindent 
Then applying the local reflexivity principle (see Lemma~3.1 of 
\cite{JRZ}), 
choose $\tilde T_\beta$ a weak* continuous finite rank operator on $X^{**}$ 
so that 
$$\|\tilde T_\beta\| < \lambda +\ep \text{ and } \tilde T_\beta |G = 
T_\alpha|G\ .
\tag 102$$
It then follows that $\tilde T_\beta \to I$ weak* on $X^{**}$, and finally so 
does the net $(\frac{\lambda}{\lambda+\ep} \tilde T_\beta)_{\beta\in\D}$.

Now choose for all $\alpha$, $S_\alpha$ a linear operator on $X^*$ with 
$S_\alpha^* = T_\alpha$, $\|S_\alpha\| \le\lambda$. 
But then it follows immediately from (ii) that 
$$S_\alpha f\to f\ \text{ weakly for all }\ f\in X^*\ .
\tag 103$$
But then there exists a net $(V_\alpha)$ of convex combinations of the 
$S_\alpha$'s so that $V_\alpha f\to f$ in norm for all $f$ in $X^*$. 
Hence $X^*$ has the $\lambda$-bap. 

(iii) $\To$ (i). 
Suppose $X^*$ has the $\lambda$-bap. 
Then it is a standard fact that $X$ has the $\lambda$-bap. 
(Actually, our proof that in Proposition~3.12, 
(ii) $\To$ (iii), already gives the argument.) 
Now let $\ep>0$, and let $F$ and $G$ be finite dimensional subspaces of $X$ 
and $X^*$ respectively. 

Then choose $S$ a finite rank operator on $X^*$ so that
$$S|F = I|F\ \text{ and }\ \|S\| <\lambda  +\ep\ .
\tag 104$$
Next, let $H = S^* X^{**}$. 
$H$ is finite dimensional so by the local reflexivity principle, we may 
choose $U:H\to X$ with 
$$\|U\| < 1+\ep \ \text{ and }\ \langle Uh,f\rangle = 
\langle h,f\rangle \text{ for all } h\in H \text{ and } f\in F\ .
\tag 105$$
Finally, let $T= US^*$. 
Then for all $g\in G$ and $f\in F$, 
$$\eqalign{
\langle Tg,f\rangle = \langle US^* g,f\rangle 
&= \langle S^* g,f\rangle\ \text{ by (105)}\cr
&= \langle g,Sf\rangle = \langle g,Sf\rangle\ \text{ by (104).}\cr}
\tag 106$$
Of course $\|T\| < (1+\ep)(\lambda +\ep)$, whence since 
$\lambda >0$ is arbitrary, $X$ is $\lambda$-ELR.\qed
\enddemo

\remark{Remark}
Extendable Local Reflexivity may easily be extended to the complete 
category, and then the quantized versions of our results are valid. 
Thus, we {\it define\/} an operator space $X$ to be Completely Extendably 
Locally Reflexive (CELR) if there is a $\lambda\ge1$ so that for all 
$\ep>0$ and finite-dimensional subspaces $F$ and $G$ of $X^*$ and 
$X^{**}$, (82) holds, except that we replace 
``$\|T\|$'' by ``$\|T\|_{\cb}$'' in (82)(iii). 
In case $\lambda$ works, we say $X$ is $\lambda$-CELR. 
We then obtain that appropriate quantized versions of Proposition~3.12 
and Theorem~3.13 are valid. 
Thus in Proposition~3.12, we replace ``Banach'' by ``operator'', 
``$\lambda$-ELR'' by ``$\lambda$-CELR'' and 
``$\|T_\alpha\|$'' by ``$\|T_\alpha\|_{\cb}$''. 
\endremark

A quantized version of Theorem 3.13 goes as follows:

\proclaim{Theorem 3.13$'$}
Let $X$ be a given operator space.
Then the following are equivalent.
\roster
\item"(i)" $X$ is {\rm CELR} and has the {\rm cbap}.
\item"(ii)" There exists a uniformly completely bounded net $(T_\alpha)$
of weak*-continuous finite-rank operators from $X^{**}$ to $X$ with
$T_\alpha x^{**}\to x^{**}$ weak* for all $x^{**}\in X^{**}$.
\item"(ii)$'$" $X$ is locally reflexive and $X$ satisfies (ii) without
assuming the $T_\alpha$'s  are weak*-continuous.
\item"(iii)" $X^*$ has the {\rm cbap} and $X$ is locally reflexive.
\endroster
\endproclaim 

Also the quantitative statements go through; if $X^*$ has the
$\lambda$-cbap and $X$ is $\beta$-locally reflexive, then $X$ is
$\lambda\beta$-CELR and has the $\lambda\beta$-cbap,
while if $X$ is $\lambda$-CELR and has the
$\beta$-cbap, then $X^*$ has the $\lambda^2\beta$-cbap.
Moreover if $(T_\alpha)$ satisfies (ii) with $\|T_\alpha\|_{\cb} \le\lambda$
for all $\alpha$, then $X$ and $X^*$ both have the $\lambda$-cbap
and $X$ is $\lambda$-CELR.
It then follows that nuclear $C^*$-algebras are 1-CELR, for it is known
that such are 1-locally reflexive
with duals having the 1 cbap which are also 1-locally reflexive [EJR].

We are indebted to N.~Ozawa for pointing out that the implication
(iii) $\Rightarrow$ (i) is false without the assumption that $X$ is
locally reflexive.
Actually, we construct a non-reflexive operator space $Y$ in
Corollary~4.9 so that $\bK_0 \subset Y\subset \bK_0^{**}$ with $Y/\bK_0$
completely isometric to $c_0$.
As pointed out to us by N.~Ozawa, since $\bK_0$ is a complete $M$-ideal
$\bK_0^{**}$, $Y^*$ is completely isometric to $C_1\oplus \ell^1$
($\ell^1$-direct sum), whence $Y^*$ has the cmap and moreover
$Y^{**}$ is isometrically injective.
It can also be seen (using arguments similar to those for Theorem~4.7
below), that $Y$ fails the cbap, thus answering a question of Ozawa's.

The next result yields an unusual connection between Extendable Local 
Reflexivity and the CSCP.

\proclaim{Theorem 3.14} 
Let $X\subset Y$ be separable operator spaces so that $X$ has the 
{\rm CSCP} and $X^{**}$ is isomorphically mixed injective. 
Suppose there exists an operator space $Z$ which is (Banach) 
{\rm ELR} and $X\subset Y\subset Z$. 
Then $X$ is complemented  in $Y$.
\endproclaim

First, an immediate consequence. 

\proclaim{Corollary 3.15} 
Suppose $B(H)$ or every separable $C^*$ algebra is {\rm ELR}, and let $X$ 
be an operator space with the {\rm CSCP} so that $X^{**}$ is mixed injective.
Then $X$ has the {\rm MSEP}; in particular, $K$ has the {\rm MSEP}.
\endproclaim

\demo{Proof}
Let $X$ have the CSCP and suppose $X\subset Y\subset B(H)$ with $Y$ separable. 
If $B(H)$ is ELR, $X$ is complemented in $Y$ by Theorem~3.14. 
Letting $\A$ be the $C^*$-algebra generated by $Y$, $\A$ is separable, so 
again $X$ is complemented in $\A$ and hence in $Y$ by Theorem~3.14.\qed
\enddemo

Theorem 3.14 is a simple consequence (via known results) of the crucial

\proclaim{Lemma 3.16} 
Let $X\subset Y$ be operator spaces. 
Assume the following:
\roster
\item"(i)" $X$ is locally reflexive
\item"(ii)" $X^{**}$ is complemented in $Y^{**}$
\item"(iii)" $Y$ is {\rm ELR}. 
\endroster
Then there is a {\rm NEW} operator space structure on $Y$, agreeing 
(isometrically) with the given one on $X$, so that $(Y,\NEW)$ 
is locally reflexive.
\endproclaim

\remark{Remark}
Our proof yields that if $X$ is $\lambda$-locally reflexive, $X^{**}$ is 
$\beta$-cocomplemented in $Y^{**}$, and $Y$ is $\gamma$-ELR, then 
$(Y,\NEW)$ is $(\gamma\beta +\lambda\beta+\lambda)$-locally reflexive.
\endremark

\demo{Proof of Theorem 3.14}
Since $X^{**}$ is mixed injective, $X^{**}$ is complemented in $Z^{**}$. 
By the Lemma, choose a NEW operator structure on $Z$ which agrees with 
that on $X$ so that $(Z,\NEW)$ is locally reflexive. 
But then $(Y,\NEW)$ is also locally reflexive, (see \cite{ER}). 
Hence $X$ is completely complemented in $(Y,\NEW)$, which of course gives 
that $X$ is complemented in $Y$.\qed 
\enddemo

We now proceed with the proof of Lemma 3.16. 
The idea goes as follows. 
By a standard Banach space construction (which we give), $X^{\bot\bot}$ 
is in fact weak*-complemented in $X^{**}$. 
Now letting $Z$ be a weak* complement, $(Y,\NEW)$ is defined in such 
a way that $(Y^{**},\NEW)$ {\it coincides\/} on $X^{\bot\bot}$ with 
the given operator space structure, while it is {\it equivalent\/} to 
$\MAX$ on $Z$.  
The hypothesis that $Y$ is ELR then allows us to obtain a ``local reflexivity 
operator'' $T_1 :G_1\to Y$, for given $G_1\subset Z$ finite-dimensional,  
with $T_1$ uniformly completely bounded, and also given $G_2\subset X^{**}$, 
$T_2 :G_2\to X$ is found by the local reflexivity of $X$. 
Then if $G= G_1\oplus G_2$, $T= T_1\oplus T_2$ is the desired local 
reflexivity operator. 

\demo{Proof of Lemma 3.16}
We identify $X^{**}$ with $X^{\bot\bot}$, $X^\bot$ with $(Y/X)^*$, and, 
as usual, $Z$ with its canonical embedding in $Z^{**}$, (for any Banach 
space $Z$). 
We first have (the standard result) that the hypotheses are equivalent to: 
$X^{**}$ is weak*-complemented in $Y^{**}$. 
In fact, fix $\beta\ge1$.
\enddemo

\definition{Fact}
$X^{**}$ is $\beta$-cocomplemented in $Y^{**}$ iff $X^\bot$ is 
$\beta$-complemented in $Y^{*}$.
\enddefinition

\demo{Proof}
Let $L:Y^{***}\to Y^*$ be the canonical projection defined by 
$$\langle Ly^{***},y\rangle = \langle y^{***},y\rangle\ \text{ for all }\ 
y^{***}\in Y^{***},\ y\in Y\ .
\tag 107$$
Now suppose first that $P:Y^{**}\to Y^{**}$ is a projection with 
$\ker P= X^{\bot\bot}$ and $\|P\| \le \beta$. 
Define $Q$ by 
$$Q=L\circ (P^*|Y^*)\ .
\tag 108$$
Now we claim that $Q$ is a projection on $Y^*$, onto $X^\bot$; of course 
it's trivial that $\|Q\| \le\beta$. 

By definition, we have for all $y^*\in Y^*$ and $y\in Y$, that 
$$\langle Qy^*,y\rangle = \langle LP^* y^*,y\rangle 
= \langle P^* y^*,y\rangle = \langle y^*,Py\rangle\ .
\tag 109$$
Now suppose first $y^*\in X^\bot$. 
Then for $y\in Y$, 
$$\eqalign{
\langle Qy^*,y\rangle & = \langle y^*,Py\rangle\ \text{ by (109),}\cr
&= \langle y^*,y\rangle \ \text{ since }\ y-Py\in X^{\bot\bot}\ .\cr}
\tag 110$$
Hence $Qy^* = y^*$. 
On the other hand, if $y^*$ is arbitrary and $x\in X$, then 
$$\langle Qy^*,x\rangle = \langle y^*,Px\rangle =0
\tag 111$$
since $x\in X^{\bot\bot}$, (109) and (111) prove our claim. 
Of course conversely if $Q:Y^*\to X^\bot$ is a projection with 
$\|Q\|\le\beta$, $Q^* \defeq P$ is a projection on $X^{**}$ with kernel 
$X^{\bot\bot}$.

Next, for $Z$ an arbitrary operator space, let ${\|\cdot \|}_{\op(Z)}$
denote the given norm on $\bK \otimes Z$; thus also 
$\|\cdot\|_{\op(Z^*)}$ is then the induced norm on $\bK\otimes Z^*$, 
given by the expression
$$\|T\|_{\op(Z^*)} = \sup \{ \|\langle T,S\rangle\| :S\in \bK\otimes Z, 
\|S\|_{\op(Z)} \le 1\}
\tag 112$$
where for $T= \sum K_i \otimes z_i^*$ in $\bK\otimes Z^*$, 
$S= \sum L_j\otimes z_j$ in $K\otimes z$, 
$$\langle T,S\rangle = \sum_{i,j} z_i^* (z_j) K_i\otimes L_j 
\tag 113$$
(regarded as an operator on $\ell^2 \otimes \ell^2$). 
Recall also, for $T$ as above, 
$$\|T\|_{\MIN} = \sup \biggl\{ \Big\| \sum_i z_i^* (z) K_i\Big\| 
: z\in \Ba (Z)\biggr\}\ .
\tag 114$$
\enddemo

Now we first define a new operator structure on $K\otimes Y^*$, and 
then let $\|\cdot\|_{\NEW}$ on $K\otimes Y$ be the one induced by this.

\definition{Definition}
For $T=\sum K_i\otimes y_i^*$ in $\bK \otimes Y^*$, set 
$$\|T\|_{N^*} = \max \left\{ \|T\|_{\MIN},\sup \{ \|\langle T,S\rangle\|: 
S\in \bK\otimes X,\ \|S\|_{\op(X)} \le 1\right\}\ .
\tag 115$$
Now it is easily verified that $\|\cdot\|_{N^*}$ on $\bK\otimes Y^{**}$ 
satisfies Ruan's axioms (cf. \cite{ER}), here $Y^*$ is 
indeed an operator space in this new structure. 
Next, we observe that $\|\cdot\|_{N^*}$ is induced by a NEW operator 
structure  on $Y$. 
It suffices to prove that given $n$, $T= \sum_{i=1}^n K_i\otimes y_i^*$, 
and a net $(T_\alpha)$ with 
$T_\alpha = \sum_{i=1}^n K_i \otimes y_{i,\alpha}^*$ 
for all $\alpha$, then if $y_{i,\alpha}^* \to y_i^*$ $\omega^*$ for all $i$ and 
$\|T_\alpha\|_{N^*} \le1$ for all $\alpha$, also $\|T\|_{N^*} \le1$. 
But it is evident that then given $y\in Y$, $\|y\|\le 1$, that 
$$\sum_i y_{i,\alpha}^* (y) K_i \to \sum y_i^* (y)K_i\ \text{ in norm,}
\tag 116$$
and moreover given $S\in \Ba(K\otimes X)$, that 
$$\langle T_\alpha,S\rangle \to \langle T,S\rangle\ \text{ in norm.}
\tag 117$$
Hence $\|\sum y_i^* (y) K_i\| \le1$, so $\|T\|_{\MIN}\le 1$, and also 
$\|\langle T,S\rangle\| \le1$, thus $\|T\|_{N^*} \le1$ as desired. 

Now let $\|\cdot\|_{\NEW}$ be the operator space structure induced on 
$K\otimes Y$ by $\|\cdot\|_{N^*}$; we have thus by duality that 
$$\|\cdot\|_{\NEW^*} = \|\cdot\|_{N^*}\ .
\tag 118$$

Next, we show that $\|\cdot\|_{\NEW}$ equals $\|\cdot\|_{\op(X)}$ on 
$K\otimes X$. 
Now first note that 
$$\|\cdot\|_{\NEW}\ge \|\cdot\|_{\op(Y)}\ .
\tag 119$$
Indeed, (119) follows immediately by duality, since 
$\|T\|_{N^*} \le \|T\|_{\op(Y^*)}$ for all $T\in\bK\otimes Y^*$. 
but if $S\in K\otimes Y$ and $\|S\|_{\op(X)} \le1$, then by definition, 
$\|\langle T,S\rangle \le1$ for all $T\in \bK\otimes Y^*$ with 
$\|T\|_{N^*} \le1$, hence $\|S\|_{\NEW} \le1$; proving 
$$\|\cdot\|_{\NEW} \le \|\cdot\|_{\op (X)}\ \text{ on }\ K\otimes Y\ .
\tag 120$$

We now assume that $X^{**}$ is $\beta$-cocomplemented in $Y^{**}$, i.e., 
$X^\bot$ is $\beta$-complemented in $Y^*$, by the Fact. 
Now choose $P:Y^*\to X^\bot$ a projection with $\|P\|\le\beta$, and let 
$E= \ker P$. 
We next claim that $\|\cdot\|_{\NEW^{**}}$ is {\it equivalent\/} to 
$\|\cdot\|_{\MAX}$ on $K\otimes E^\bot$. 
Now it follows immediately from the definition that 
$$\|T\|_{\NEW^*} = \|T\|_{\MIN}\ \text{ for all }\ T\in K\otimes X^\bot\ .
\tag 121$$
By duality, we have that for any $S\in K\otimes Y^{**}$, 
$$\|S\|_{\MAX} = \sup \{ \|\langle T,S\rangle\| : T\in K\otimes Y^*,\ 
\|T\|_{\MIN} \le 1\}\ .
\tag 122$$
But if $T\in K\otimes Y^*$, say $T= \sum K_i\otimes y_i^*$, then 
letting $\tilde PT = \sum K_i\otimes Py_i^*$, we have that 
$$\|\tilde PT\|_{\MIN} \le \beta \|T\|_{\MIN}
\tag 123$$
and of course  $\tilde P T\in K\otimes X^\bot$. 
Hence we obtain that if $S\in K\otimes E^\bot$, then for any 
$T\in K\otimes Y^*$ 
$$\langle S,T\rangle  = \langle S,\tilde PT\rangle 
\tag 124$$
whence 
$$\eqalign{
\|\langle S,T\rangle\| = \|\langle S,\tilde PT\rangle \|
&\le \|S\|_{\NEW^{**}} \|\tilde P T\|_{\NEW^*}\cr
&= \|S\|_{\NEW^{**}} \|\tilde P T\|_{\MIN}\ .\cr}
\tag 125$$
Thus 
$$\|S\|_{\MAX} \le \beta \|S\|_{\NEW^*} 
\tag 126$$ 
as desired. 
Finally, we show that $(Y,\NEW)$ is locally reflexive. 
Assume then that $X$ is $\lambda$-locally reflexive and now suppose $Y$ 
is $\gamma$-ELR. 
Let $F,G$ be finite-dimensional spaces with $F\subset Y^*$, $G\subset Y^{**}$. 
Now we may assume without loss of generality, by simply enlarging $G$ 
and $F$ if necessary, that 
$$G = G_1 \oplus G_2\ \text{ and }\ F=F_1\oplus F_2
\tag 127$$
with 
$$G_1 \subset E^\bot,\ G_2\subset X^{\bot\bot},\ F_1\subset X^\bot,\ 
F_2\subset E\ .
\tag 128$$
Let $\ep >0$. 
Since $X$ is $\lambda$-locally reflexive, choose $T_2 :G_2\to X$ with 
$$\eqalign{
\text{(i)}\hfil\quad&\|T_2\|_{\cb} <\lambda +\ep\cr
\text{(ii)}\hfil\quad&\langle T_2y,f\rangle = \langle g,f\rangle
\text{ for all } g\in G_2,\ f\in E\ .\cr}
\tag 129$$
Since $Y$ is $\lambda$-ELR, we may choose $T_1:E^\bot \to Y^{**}$ with 
$$\eqalign{
\text{(i)}\hfil\quad&\|T_1\| < \gamma+\ep\cr
\text{(ii)}\hfil\quad&T_1G_1\subset Y\cr
\text{(iii)}\hfil\quad&\langle T_1g,f\rangle = \langle g,f\rangle 
\ \text{ for all }\ g \in G_1,\ f\in F_1\ .\cr}
\tag 130$$
Now it follows that 
$$\| T_1\|_{\cb} < (\gamma +\ep) \beta\ .
\tag 131$$
(Here, we are computing the cb norm with respect to  $(E^\bot,\NEW^{**})$.) 
Indeed, if $I$ denotes $I|K$, then if $S\in K\otimes E^\bot$, 
$$\|(I\otimes T_1)(S)\| \le \|T\|\, \|S\|_{\MAX} 
< (\gamma +\ep) \beta \|S\|_{\NEW^{**}}\ . 
\tag 132$$
Finally, we define  $T:G \to Y$ by 
$$T= T_1|G_1\oplus T_2\ .
\tag 133$$
Then by (129)(ii), (130)(iii), and (127), 
$$\langle Tg,f\rangle = \langle g,f\rangle\ \text{ for all }\ 
g\in G,\ f\in F\ .
\tag 134$$

Now letting $R= P^*|G$, we have 
$$\|R\|_{\cb} \le \|P^*\|_{\cb} = \|P\|_{\cb} 
\le \|P\| \le \beta 
\tag 135$$
since $X^\bot$ has the MIN operator structure by (121). 
Now $T= T_1 R + T_2 (I-R)$ by (128), hence 
$$\eqalign{
\|T\|_{\cb} &\le \|T_1R\|_{\cb} + \|T_2 (I-R)\|_{\cb}\cr 
&\le (\gamma +\ep) \beta + (\lambda +\ep)(1+\beta)\cr}
\tag 136$$
by (129), (131) and (135). 

Since $TG\subset Y$ and (122) holds, we have established that 
$(Y,\NEW)$ is $(\gamma+\lambda)\beta+\lambda$-locally reflexive.\qed
\enddefinition

\remark{Remarks}
1. The alert reader may notice that the ELR assumption on $Y$ is used only at 
the very end. 
Thus, without this, we still obtain that $(Y,\NEW)$ coincides on $X$ with the 
original operator space structure, and $(Y^*,\NEW)$ is still the MAX 
structure on $E^\bot$ (to a constant), and the given structure on $X^{**}$. 
However if $G\subset E^\bot$ is finite-dimensional, we cannot insure that 
a Banach local reflexivity operator $T:G\to Y$ is uniformly completely 
bounded, since $G$ may not have MAX as its induced operator structure. 
The synthesis of the ELR concept occurred precisely to overcome this 
(apparently insurmountable) difficulty. 

2. There is really no reason to assume that $Y$ is an operator space at all. 
We really make no essential use of the given operator space structure 
on $Y$; the inequality (110) can instead be easily established directly 
(replacing "$X$'' in its statement). 
We also obtain that $X^{**}$ is completely complemented in $(Y^{**},\NEW)$ 
(in fact if $X^{**}$ is $\beta$-cocomplemented in $Y^{**}$, it is completely 
$\beta$-cocomplemented in $(Y^{**},\NEW)$). 

3. It is an open question 
{\it if maximal operator spaces are locally reflexive\/}. 
If the answer to this question is affirmative, the conclusion of 
Lemma~3.16 would hold {\it without\/} the assumption that $Y$ is ELR; 
consequently Theorem~3.14 would hold {\it without\/} the assumption 
of the existence of the ELR $Z$ in its statement, and it would follow 
that $\bK$ has the MSEP. 
(Moreover here, we would just require that separable maximal operator 
spaces are locally reflexive.) 
Indeed, the NEW operator space structure on $Y$ is defined so that the 
induced structure on $E^\bot \subset Y^{**}$ is equivalent to MAX there. 
ELR of $Y$ is used solely to insure the existence of the ``local 
reflexivity'' of $T_1|G_1$ with controlled cb-norm. 
Now suppose $(Y,\MAX)$ is locally reflexive. 
But then we could choose $T_1:G_1\to Y$ satisfying (130)(iii) with 
$\|T_1\|_{\cb}\le \tau$, where $\tau$ is a constant depending only on the 
local reflexivity constant of $(Y,\MAX)$ and on $\beta$ (as defined in the 
proof). 
We note concerning this open question that it is equivalent (in general 
to the problem: 
{\it is $(B(H),\MAX)$ locally reflexive\/}? 
Indeed, fixing a maximal operator space $Y$, choose a Hilbert space $H$ 
so that $Y\subset B(H)$. 
But {\it then\/} the induced operator structure on $Y$ via $(B(H),\MAX)$ 
coincides with the given maximal structure, thanks to the injectivity 
of $B(H)$. 
Thus if $(B(H),\MAX)$ is locally reflexive, so is $Y$. 
\endremark
\bigskip

\centerline{{\smc Section} 4}
\centerline{\bf $\bK_0$ fails the CSEP: a new proof and generalizations}
\smallskip

To formulate the main result of this section, we first recall a concept 
introduced in \cite{R2}. 

\definition{Definition}
A family $\Z$ of operator spaces is said to be {\it finite matrix type\/} 
if there is a $C\ge 1$ so that for any finite-dimensional operator space 
$E$, there is  an $n=\bn(E)$ so that  
$$\|T\|_{\cb} \le C\|T\|_n \text{ for all linear operators } 
T:E\to Z\text{ and all } Z\in\Z\ .$$
If $C$ works, we say that $\Z$ is of $C$-finite matrix type, or briefly, that 
$\Z$ is $C$-finite. 
Finally, we say that an operator space $Z$ is of finite matrix type provided 
$\{Z\}$ has this property. 
\enddefinition

It is established in \cite{Ro2}, Proposition 2.15, that 
{\it if for some $\lambda$, $Z_1,Z_2,\ldots$ are separable 
$\lambda$-injective operator spaces with $\{Z_1,Z_2,\ldots\}$ of finite 
matrix type, then $(Z_1\oplus Z_2\oplus \cdots)_{c_0}$ has the\/} CSEP. 
Our main result in this section establishes the converse. 

\proclaim{Theorem 4.1} 
Let $Z_1,Z_2,\ldots$ be operator spaces so that $\{Z_1,Z_2,\ldots\}$ is not 
of finite matrix type. 
If all of the $Z_i$'s have finite matrix type, let $Z= (Z_1\oplus Z_2
\oplus \cdots)_{c_0}$. 
Otherwise, choose $i$ so that $Z_i$ is not of finite matrix type, and set 
$Z= c_0(Z_i)$. 
Then there exists an operator space $Y$ with $Y/Z$ separable such that 
$Z$ is not completely complemented in $Y$.
\endproclaim

We then easily obtain a converse to the result from \cite{Ro2} mentioned 
above, in view of the fact that separable injective operator spaces are 
necessarily injective (Corollary~2.9 of \cite{Ro2}). 

\proclaim{Corollary 4.2} 
Let $Z_1,Z_2,\ldots$ be reflexive separable operator spaces so that 
$(Z_1\oplus Z_2\oplus \cdots)_{c_0}$ has the {\rm CSEP}, and assume that 
$Z_i$ is of finite matrix type for all $i$. 
Then there is a $\lambda$ so that $Z_j$ is $\lambda$-injective for all $j$, 
and $\{Z_1,Z_2,\ldots\}$ is of finite matrix type.
\endproclaim

\remark{Remark} 
We conjecture that the last hypothesis is superfluous; see the Conjecture 
following Corollary~4.3.
\endremark

\demo{Proof of 4.2} 
By the results of \cite{Ro2}, there exists a $\lambda$ so that 
$(Z_1\oplus Z_2\oplus \cdots)_{c_0}$ has the $\lambda$-CSEP. 
Hence for each $j$, $Z_j$ has the $\lambda$-CSEP. 

Since $Z_j$ is reflexive, $Z_j$ is $\lambda$-injective by Proposition~2.10 
of \cite{Ro2}. 
Of course Theorem~4.1 then yields that $\{Z_1,Z_2,\ldots\}$ is of 
finite matrix type.\qed
\enddemo

Now standard results yield that $\{\M_1,\M_2,\ldots\}$ is {\it not\/} of finite 
matrix type, where for all $n$, $\M_n$ denotes the operator space of 
$n\times n$ matrices.  
(We give a quantitative refinement of this fact below.) 
Thus we obtain the result of E.~Kirchberg \cite{Ki1} (see also \cite{W}):

\proclaim{Corollary 4.3} 
$\K_0$ fails the {\rm CSEP}.
\endproclaim

\example{Conjecture} 
If a separable operator space has the CSEP, it is of finite matrix type.
\endexample

The next immediate consequence of 4.1 supports this conjecture.

\proclaim{Corollary 4.4} 
Let $Z$ be a separable operator space which is not of finite matrix type. 
Then $c_0(Z)$ fails the {\rm CSEP}. 
Hence if $c_0(Z)$ is completely isomorphic to $Z$, then $Z$ fails the 
{\rm CSEP}.
\endproclaim

We now proceed with the proof of Theorem 4.1. 
The following construction gives the crucial tool. 

\proclaim{Lemma 4.5} 
Let $(Z_1,Z_2,\ldots)$ be a given sequence of operator spaces, $k$ a positive 
integer, $C>1$, and $E$ an $m$-dimensional operator space. 
Assume there exists a sequence $1=n_0<n_1<n_2<\cdots$ of positive integers 
and for all $k\ge1$, a linear map $U_k :E\to Z_k$ so that 
\roster
\item"(137i)\ " \hskip1.5truein $\|U_k\|_{\cb} \le1$
\item"(137ii)$\,$" \hskip1.5truein $\|U_k\|_{n_k} > 1- \frac1k$
\item"(137iii)" \hskip1.35truein $\|U_k\|_{n_{k-1}} \le \frac1C\ \frac{k}{k-1}$ if $k>1$.
\endroster
Then setting $Z= (Z_1\oplus Z_2\oplus \cdots)_{c_0}$, there exists an 
operator space $Y\supset Z$ with $\dim Y/Z \le m$ so that $\|P\|_{\cb}\ge C$ 
for any surjective linear projection $P:Y\to Z$.
\endproclaim 

\demo{Proof}
In this discussion, we let $I_j$ denote the identity map on $\M_j$. 
We construct $Y$ as a subspace of 
$W\defeq (Z_1\oplus Z_2\oplus \cdots)_\infty$.
Define $U:E\to W$ by 
$$U(x) = (U_1(x),U_2(x),\ldots)\ \text{ for }\ x\in E
\tag 138$$
and let $F= U(E)$, $Y= Z+F$. 
Let $P:Y\to W$ be a linear projection. 
Now given $k_0>1$, by making a small perturbation if necessary, we may 
without loss of generality assume that there is a $k>k_0$ with 
$$P\!(F)\subset Z_1\oplus\cdots \oplus Z_{k-1}\ .
\tag 139$$
Let $Q_j$ be the coordinate projection from $W$ onto $Z_j$, for all $j$. 
Now by (137ii), choose $\tau \in E\otimes \M_{n_k}$ with 
$\|\tau \|=1$ and 
$$\|U_k\otimes I_{n_k} (\tau)\| > 1-\frac1k\ .
\tag 140$$
Then letting $\beta = (U\otimes I_{n_k})(\tau)$, we have by (137i) and (140) 
that 
$$1\ge \|\beta\| \ge \|Q_k \otimes I_{n_k}(\beta)\| > 1-\frac1k\ ,
\tag 141$$ 
and by (137iii) that 
$$\|Q_\ell \otimes I_{n_k} (\beta)\| \le \frac1C \ \frac{k}{k-1} \text{ for 
all } \ell > k\ .
\tag 142$$ 

Finally, let $\gamma =\beta - \sum_{j=1}^k (Q_j \otimes I_{n_k})(\beta)$. 
Then 
$$\|\gamma\| = \sup_{j>k} \|Q_j\otimes I_{n_k}(\beta)\| 
\le \frac1C\ \frac{k}{k-1} \le \frac1C\ \frac{k_0}{k_0-1}\ . 
\tag 143$$

However we have that 
$(Q_k\otimes I_{n_k})(P\otimes I_{n_k})(\gamma) = -Q_k\otimes I_{n_k}(\beta)$ 
by (134), hence 
$$\align 
\|P\otimes I_{n_k} (\gamma)\| 
& \ge \|(Q_k \otimes I_{n_k}) (P\otimes I_{n_k})(\gamma)\| 
\tag 144\\
&= \|Q_k\otimes I_{n_k}(\beta)\| > 1-\frac1k \text{ by (141)} \\
&\hskip1.3truein \ge 1-\frac1{k_0}\ .
\endalign$$
Since $k_0>1$ is arbitrary, (143) and (144) yield that 
$\|P\|_{\cb} \ge C$, as desired.\qed
\enddemo

The next quantitative result easily yields Theorem 4.1.

\proclaim{Lemma 4.6} 
Let $C>1$ and let $\Z$  be a family of operator spaces which is not 
$C$-finite. 
There exist $Z_1,Z_2,\ldots$ in $\Z$ and an operator space $Y\supset Z
\defeq (Z_1\oplus Z_2\oplus \cdots)_{c_0}$ with $Y/Z$ finite-dimensional 
so that $\|P\|_{\cb} \ge C$ for any linear surjective projection $P:Y\to Z$. 
\endproclaim

\demo{Proof} 
Choose $E$ a finite dimensional operator space so that for all $n\in \nat$, 
there exists a $Z\in \Z$ and a linear operator $U:E\to Z$ with 
$$\|U\|_{\cb} =1\quad\text{and}\quad \|U\|_n < \frac1C\ .
\tag 145$$

Also note, that for any completely bounded map $T$ between operator spaces, 
$$\|T\|_{\cb} = \sup_n \|T\|_n \ .
\tag 146$$ 

Using (145), choose $Z_1\in\Z$ and a linear operator $U_1: E\to Z_1$ 
with $\|U_1\|_{\cb} =1$ and $\|U_1\|_1<\frac1C$. 
Choose $n_1>1$ with $\|U_1\|_{n_1} >0$. 
Suppose $k>1$ and $n_{k-1}$ has been chosen. 
By (145), we may choose $Z_k \in\Z$ and a linear operator 
$U_k :E\to Z_k$ with $\|U_k\|_{\cb}=1$ and $\|U_k\|_{n_{k-1}} <\frac1C$. 
Then using (146), choose $n_k >n_{k-1}$ with $\|U_k\|_{n_k} >1-\frac1k$. 

This completes the inductive construction. 
Then (137) holds for all $k$, so the $U_k$'s satisfy the hypotheses of 4.5, 
which thus yields Lemma~4.6.\qed
\enddemo 

We are now prepared for the 

\demo{Proof of Theorem 4.1}
Suppose first that $X$ is an operator space which is not of finite matrix 
type, and let $Z = c_0(X)$. 
Then by Lemma~4.6, for each $n\in\nat$ we may choose $Y_n$ an operator space 
with $Y_n\supset Z$ so that $Z$ is not $n$-completely complemented in $Y_n$ 
and $Y_n/Z$ is finite dimensional. 
Let $Y= (Y_1\oplus Y_2\oplus\cdots)_{c_0}$ and $\tilde Z = (Z\oplus Z\oplus 
\cdots)_{c_0}$. 
Then $\tilde Z$ is canonically completely isometric to $Z$, $Y/\tilde Z$ 
is separable, and $\tilde Z$ is not completely complemented in $Y$. 

Now let $\{Z_1,Z_2,\ldots\}$ be as in the statement of 4.1. 
If $i$ is such that $Z_i$ is not of finite matrix type, the above argument 
establishes the conclusion of 4.1. 
Otherwise, Lemma~4.6 and its proof yield that we may choose infinite 
pairwise disjoint subsets $M_1,M_2,\ldots$ of $\nat$ so that for each $j$, 
letting $W_j = (\bigoplus_{i\in M_j} Z_i)_{c_0}$, there exists an operator 
space $Y_j \supset W_j$ with $Y_j/W_j$ finite-dimensional and $W_j$ not 
$j$-completely complemented in $Y_j$. 

Indeed, it follows from the  definition of families of finite matrix type 
that there then exist $\ell_1<\ell_2<\cdots$ so that $Z_{\ell_j}$ is 
{\it not\/} $j$-finite for all $j$. 
Then let $\tilde M_1,\tilde M_2,\ldots$ be infinite 
pairwise disjoint sets so that 
$\bigcup_{j=1}^\infty \tilde M_j = \{\ell_1,\ell_2,\ldots\}$. 
Now it follows that letting $\Z_j = \{Z_m :m\in\tilde M_j\}$, then $\Z_j$ 
is not of finite matrix type for all $j$; now Lemma~4.6 yields 
an appropriate infinite $M_j\subset \tilde M_j$, for all $j$, satisfying 
the above. 

Now letting $\tilde Y = (Y_1\oplus Y_2\oplus\cdots)_{c_0}$ and 
$\tilde W= (W_1\oplus W_2 \oplus \cdots)_{c_0}$, then $\tilde Y/\tilde W$ 
is separable and $\tilde W$ is not completely complemented in $Y$. 
Finally, let $M_0 = W\sim \bigcup_{j=1}^\infty M_j$ and 
$\tilde Z = (\bigoplus_{i\in M_0} Z_i)_{c_0}$. 
Then let $Y = \tilde Y\oplus \tilde Z$. 
$Z$ is canonically isometric to $\tilde W\oplus \tilde Z$, and of 
course $\tilde W\oplus \tilde Z$ is uncomplement in $Y$ and 
$Y/(\tilde W\oplus\tilde Z) = \tilde Y/\tilde W$.\qed
\enddemo

We now give a ``tight'' quantitative version of Corollary~4.3 (which is one 
of the main motivating results of this section). 
Recall that for a finite-dimensional operator space $X$, the exactness 
constant of $X$, denoted $\Ex(X)$, is defined by 
$$\align 
\Ex (X) & = \inf \{d_{\cb} (X,F) : F\subset \bK\}\\
& \faceq \inf  \{d_{\cb} (X,F) : F\subset \M_n\text{ for some } n\}\ .
\endalign$$

\proclaim{Theorem 4.7} 
Let $E$ be a finite-dimensional operator space, and let $C= \Ex(E^*)$. 
There exists an operator space $Y$ containing $\bK_0$ and a 
finite-dimensional subspace $F$ of $Y$ so that 
\roster
\item"(i)" $Y/\bK_0$ is completely isometric to $E$.
\item"(ii)" $\bK_0$ is Banach $(1+\ep)$ co-completely in $Y$ for 
every $\ep >0$.
\item"(iii)" $\|P\|_{\cb} \ge C$ for any surjective linear projection 
$P:Y\to \bK_0$.
\endroster
\endproclaim 

We first require a lemma, which really yields a precise local, 
quantitative version of the fact that $\{\M_n :n=1,2,\ldots\}$ is {\it not\/} 
of finite matrix type.

\proclaim{Lemma 4.8} 
Let $E$ be a finite-dimensional operator space, $\ell>1$, $\ep>0$ and 
set $C=\Ex (E^*)$. 
There exist an $m$ and a {\rm 1--1} operator $T:E\to \M_{m}$ satisfying 
the following:
\roster
\item"(i)" $(1+\ep) C > \|T\|_{\cb} > (1-\ep) C$
\item"(ii)" $\frac1{1+\ep} \|x\| \le \|T\otimes I_\ell (x)\| \le 
(1+\ep) \|x\|$ for all $x\in E\otimes \M_\ell$.
\endroster
\endproclaim 

\demo{Proof} 
We let $P_k:\M_\infty \to \M_k \subset\M_\infty$ be the natural 
truncation operator; i.e., 
\roster
\item"{}" $P_k(a_{ij}) = a_{ij}$ if $1\le i,j\le k$
\item"{}" $P_k (a_{ij}) = 0$ otherwise.
\endroster
Of course $\|P_k\|_{\cb}=1$ and $P_kT\to T$ in the strong operator topology.

We first note that it suffices to find $T:E\to \M_\infty$ satisfying (i) 
and (ii). 
Indeed, if such a $T$ has these properties, then for $m$ large enough, 
(since $E$ is finite-dimensional), 
$(1+\ep) C > \|P_mT\|_{\cb} > (1-\ep) C$ also and 
$$\frac1{(1+\ep)^2} \|x\| \le \|P_m T\otimes I_\ell (x)\| \le 
(1+\ep)^2\|x\|
\tag 147$$ 
for all $x\in E\otimes \M_\ell$, 
hence $\tilde T\defeq P_m T$ has the desired property (for a little 
bigger $\ep$). 

Now we dualize; without loss of generality $E^* \subset \M_\infty$. 
Next we claim that for $k$ sufficiently large, 
$$\frac1{1+\ep} \|x\| \le \|P_k \otimes I_\ell (x)\|\text{ for all } 
x\in E^* \otimes \M_\ell
\tag 148$$
and
$$\|(P_k|E^*)^{-1}\|_{\cb} < (1+\ep) C\ .
\tag 149$$
(Note that by (148), we will have $P_k|E^*$ is 1--1; setting 
$G_k^* = P_k(E^*)$, $(P_k|E^*)^{-1}$ refers to the inverse of 
$E^* @> P_k|E^*>> G_k^*$). 
Indeed, we may choose $n\ge \ell$ and $Y\subset\M_n$ with 
$d_{\cb} (E^*,Y) < (1+\ep) C$. 
Hence we may choose $T:E^* \to Y$ a linear operator with 
$$\|T\|_{\cb} \|T^{-1}\|_{\cb} < (1+\ep)^{1/2} C\ .
\tag 150$$
Next, since $E^*$ is finite-dimensional, so is $E^*\otimes \M_\ell$, so 
we can in fact choose $k$ so that 
$$\frac1{(1+\ep)^{1/2}} \|x\| \le \|P_k\otimes I_n(x)\| \text{ for all } 
x\in E^* \otimes \M_n\ .
\tag 151$$
which gives (149) immediately. 
But then we have that $T(P_k|E^*)^{-1} :G_k \to \M_n$, thus using a 
Lemma of Roger Smith (cf.\ \cite{S}, also see \cite{Pi3}), 
$$\align 
\|T(P_k|E^*)^{-1}\|_{\cb} 
& = \|T(P_k|E^*)^{-1}\|_n \tag 152\\
&\le \|T\|_{\cb} \|(P_k|E^*)^{-1}\|_n\\
&\le (1+\ep)^{1/2} \|T\|_{\cb} \text{ (by (151)}\ .
\endalign$$
But then 
$$\align
\|(P_k|E^*)^{-1}\|_{\cb} 
& = \|T^{-1} T(P_k|E^*)^{-1}\|_{\cb}\tag 153\\
& \le \|T^{-1}\|_{\cb} \|(P_k|E^*)^{-1}\|_{\cb}\\
&\le (1+\ep) C\text{ by (151) and (152)),}
\endalign$$
proving (149). 

Finally, set $G_k^* = P_k (E^*)$, let $S= (P_k|E^*)^{-1} :G_k^* \to E^*$, 
and let $T= S^* :E\to G_k$. 
Then since $\|S^{-1}\|_{\cb} = \|P_k|E^*|_{\cb} \le1$, 
and $E^* \subset \M_n$, $\|S\|_{\cb} \ge C$; hence in fact 
$$C\le \|T\|_{\cb} 
\tag 154$$
and, by (148), 
$$\frac1{1+\ep} \|x\| \le \|T\otimes I_\ell (x)\| \le \|x\| 
\text{ for all }x\in E\otimes \M_\ell\ .
\tag 155$$ 
(Also $\|T^{-1}\|_{\cb} = \|S^{-1}\|_{\cb} = \|P_k|E^*\|_{\cb} \le 1$, but 
we don't use this.) 
Thus $T$ satisfies (i) and (ii) (regarding $T(E)\subset\M_\infty$), so 
at last we obtain the desired operator by our initial observations.\qed
\enddemo

\remark{Remark} 
Buried in this  proof, we have a rather remarkable fact: 
if $X$ is a finite-dimensional subspace of $\M_\infty$, 
{\it then for $k$ sufficiently large, $P_k|X$ is {\rm 1--1} and 
$\|P_k|X\|^{-1} \to \Ex)(X)$ as $k\to \infty$.} 
That is, not only do we locate a specific $Y\subset \M_k$ with $d_{\cb}(Y|X)$ 
close to $\Ex(X)$, we also obtain that setting $Y=P_k|X$, 
$T= P_k|X : X\to Y$ satisfies $\|T\|_{\cb} \|T^{-1}\|_{\cb}$ is almost 
equal to $\Ex (X)$. 
(This fact may also be found buried in the discussion in \cite{Pi2}.)
\endremark

We are now prepared for the 

\demo{Proof of Theorem 4.7} 
Let $0<\eta <1$ with $\frac{1+\eta}{1-\eta} <1+\ep$. 
Using Lemma 4.8, we choose $1=n_0<n_1<n_2<\cdots$ and for all $k$, 
linear maps $U_k: E\to \M_{n_k}$ as follows: 
First choose $n_1>1$ and an operator $T_1:E\to \M_1$ so that (i) and (ii)
of 4.8 hold for ``$T$''~$=T$, $\ep =\frac{\eta}2$, $\ell=1$. 

Set $U_1 = T_1/\|T_1\|_{\cb}$. 
Suppose $k>1$ and $n_{k-1}$ has been defined. 
Choose $n_k>n_{k-1}$ and an operator $T_k:E\to \M_{n_k}$ so that (i) 
and (ii) of 4.8 hold for ``$T$''~$=T_k$, $\ep = \frac{\eta}{2^k}$, 
$\ell=n_{k-1}$.
Then set $U_k = T_k/\|T_k\|_{\cb}$.

This completes the inductive construction of the $U_k$'s. 
We then have for all $k$, letting 
$\tau_k = (1+\frac{\eta}{2^k}) |(1-\frac{\eta}{2^k})$ and noting that $1-\ep<
\frac1{1+\ep}$ if $\ep<1$, that 
$$\gather
\|U_k\|_{\cb} = \|U_k\|_{n_k} =1\tag 155i)\\
\|U_k\|_{n_{k-1}} = \|U_k\otimes I_{n_{k-1}}\| \le \tau_k/C\tag 155ii\\
\|(U_k\otimes I_{n_{k-1}})^{-1}\| \le \tau_k C\ .\tag 155iii
\endgather$$

Setting $Z_k = \M_{n_k}$, $C$ and the $U_k$'s fulfill the hypotheses 
of Lemma~4.5, so let $\tilde Y$ be the space given in that construction 
and simply let $Y=\tilde Y\oplus (\M_{i_1}\oplus \M_{i_2}\oplus \cdots)_{c_0}$
where $i_1< i_2<\cdots$ is an increasing enumerator of $\nat\sim \{n_1,n_2,
\ldots\}$. 
Now it is immediate that $Y$ satisfies (iii) of 4.7; let us verify the 
other assertions of 4.7 (which immediately reduce to considering 
$\tilde Y$ instead). 

Let $U$ and $Z$ be as in the proof of 4.5. 
Let $\pi :\tilde Y \to \tilde X/Z$ be the quotient map. 
Then we have for $\ell\ge1$ and $x\in E\otimes \M_\ell$, that for all 
$k>\ell$, since then $n_{k-1} \ge \ell$, 
$$\frac1{\tau_kC}\|x\| \le \|U_k\otimes I_\ell (x)\| \le (\tau_k/C)\|x\| \ .
\tag 156$$
But then for any $w\in Z\otimes\M_\ell$, 
$$\lim_{k\to\infty} 
\|Q_k \otimes I_\ell (U\otimes I_\ell (x) -w)\| = \frac{\|x\|}C\ .
\tag 157$$
Since $\ell$ is arbitrary, this shows that 
$$\|(\pi U)\otimes I_\ell (x) \| = \frac{\|x\|}C\ .
\tag 158$$
That is, $C\pi U$ is a complete isometry, proving (i) of Theorem 4.7.

Finally, let $\ep>0$, choose $k_0$ so that $\tau_k<1+\ep$ if $k>k_0$, and 
define $V:E\to \tilde Y$ by 
$$\align 
Q_j V(e) & = 0 \text{ if } j<k_0\ ,\tag 159\\
Q_j V(e) & = QU_k (e)\text{ if } j\ge k_0\ . 
\endalign$$ 
Then setting $F= V(E)$, we have that $F\oplus \bK_0 = \tilde Y$ and 
$f\in F$ and $z\in \bK_0$ imply  
$$\align 
\|f+z\| &\ge \lim_{k\to\infty} \|Q_k (f+z)\|\tag 160\\
& = \lim_{k\to\infty} \|Q_k(f)\|\\
& \ge \frac1{1+\ep} \|f\|
\endalign$$
by (156) (for $\ell=1$), 
showing that $Z$ is $(1+\ep)$-cocomplemented in $\tilde Y$, 
completing the proof.\qed
\enddemo

We now draw some immediate consequences of Theorem 4.7 and previously 
known results. 

\proclaim{Corollary 4.9} 
{\rm (a)} For all $n$, there exists an operator space $Y_n$ containing 
$\bK_0$ so that $Y_n/\bK_0$ is completely isometric to $\ell_n^\infty$ and 
$\|P\|_{\cb} \ge \sqrt{n}/2$ for any surjective linear projection 
$P:Y_n\to \bK_0$.

{\rm (b)} There exists an operator space $Y$ containing $\bK_0$ so that 
$Y/\bK_0$ is completely isometric to $c_0$ and $\bK_0$ is completely 
uncomplemented in $Y$.
\endproclaim

\demo{Proof}
(a) Set $E= \ell_n^\infty$ in Theorem 4.7. 
Then $E^* = (\ell_n^1,\MAX)$ and it is known that $\Ex (\ell_n^1,\MAX)
\ge \sqrt{n}/2$ \cite{Pi2}. 

(b) Let $Y = (\oplus Y_n)_{c_0}$, and $\tilde\bK_0 = c_0(\bK_0)$. 
Of course $\tilde\bK_0$ is isometric to $\bK_0$, $\tilde \bK_0$ 
is completely uncomplemented in $Y$ by (a), and $Y/\tilde\bK_0$ is 
completely isometric to $(\oplus \ell_n^\infty)_{c_0}$ which is 
completely isometric to $c_0$.\qed
\enddemo

\remark{Remarks} 
1. By a standard result, there exists a linear projection $P:Y_n\to \bK_0$ 
with $\|P\|_{\cb} \le \sqrt{n} +1$. 
Thus the order of magnitude result in (a) is best possible. 
Our construction yields that $\bK_0$ is Banach $(1+\ep)$-co-complemented in 
$Y_n$ and $Y$, for any $\ep>0$.

2. Actually, in part (a), we may replace $\ell_n^\infty$ by any 
$n$-dimensional Banach space $E$ endowed with the minimal operator space 
structure. 
Then by a result of M.~Junge and G.~Pisier, $Ex(E^*,\MAX) \ge \sqrt{n}/4$
\cite{JP}. 
Hence we obtain an operator space $Y_n$ containing $\bK_0$ so that $Y_n/\bK_0$ 
is completely isometric to $E$ and $\bK_0$ is not $\lambda$-completely 
complemented in $Y_n$ if $\lambda< \sqrt{n}/4$.

3. A separable operator space $X$ is defined to be {\sl nuclear\/} if 
there exists a sequence $(T_n)$ of finite rank operators on $X$ with 
$T_n\to I_X$ in the strong operator topology, so that for all $n$, there 
exist $\ell_n$ and complete contractions $U_n:X\to \M_{\ell_n}$ and 
$V_n:\M_{\ell_n}\to X$ with $T_n = V_nU_n$. 
Thus, a separable $C^*$-algebra is nuclear precisely when it is a nuclear 
operator space. 
It follows from the results of E.~Kirchberg in \cite{Ki2} that 
{\it the space $Y$ in {\rm (b)} is not nuclear; however 
$\bK_0$ and $Y/\bK_0$ are obviously nuclear\/}. 
This is in marked contrast with the algebraic case (in fact since 
$\bK_0$ {\it is\/} already an ideal in $\bK_0^{**}$, if $\bK_0\subset \A
\subset \bK_0^{**}$ with $\A$ a $C^*$-algebra, then $\A/\bK_0$ nuclear 
{\it implies\/} $\A$ is nuclear). 
Indeed, the work in \cite{Ki2} yields that were $Y$ nuclear, $Y$ would 
be 1-locally reflexive, whence $\bK_0$ would be completely complemented 
in $Y$ since $\bK_0$ has the CSCP (\cite{Ro2}), contradicting 
Corollary~4.9(b). 
\endremark

\proclaim{Corollary 4.10} 
Let $Y_n$ be as in part (a). 
Then $Y_n$ is not $\lambda$-locally reflexive for $\lambda\le (\sqrt{n}/2)-3$.
\endproclaim

\remark{Remark}
Of course $Y_n$ is locally reflexive; 
in fact just because $\dim Y_n/\bK_0 =n$, there is an absolute 
constant $c$ so that $Y_n$ is $c\sqrt n$ locally reflexive. 
\endremark

\demo{Proof} 
Suppose that $Y_n$ is $C$-locally reflexive. 
By Sublemma~3.11 of \cite{Ro2}, since $K_0^* = B(H)$ is 
isometrically injective, $\bK_0$ is $C+3+\ep$-completely complemented 
in $Y_n$ for all $\ep>0$. 
Hence by Corollary~4.9, $C+3+\ep \ge \sqrt{n}/2$ for all such $\ep>0$, 
so $C\ge (\sqrt{n}/2)-3$.\qed
\enddemo

Our next (and final) application of the arguments for Theorem~4.7 yields 
that every descending sequence of 1-exact Banach isometric 
finite-dimensional spaces is bounded below. 

\proclaim{Proposition 4.11} 
Let $(\lambda_k)$ be a sequence of real numbers with $\lambda_k\ge1$ 
for all $k$ and $\prod_{k=1}^\infty \lambda_k <\infty$. 
Let $(E_j)$ be a sequence of $1$-exact finite dimensional 
operator spaces so that 
$E_k$ is $\lambda_k$-semi-isometric to $E_{k+1}$ for all $k$. 
Then $\varlimsup_{k,n\to\infty} d_{\cb}(E_k,E_n) \le 4$. 
\endproclaim

\demo{Proof} 
For each $k$, choose $J_k :E_k\to E_{k+1}$ a linear map with 
$$\|J_k^{-1}\| \le \lambda_k\ \text{ and }\ 
\|J_k\|_{\cb}=1\ \text{ for all }\ k\ .
\tag 161$$
Suppose the conclusion were false; then by passing to a subsequence if 
necessary, we may assume that for some $C > 4$, 
$$\|J_k^{-1}\|_{\cb} \ge C\ \text{ for all }\ k\ .
\tag 162$$
(Note that if $n_1<n_2<\cdots$ is given, then letting $\tilde J_k = 
(J_{n_{k+1}-1})\cdots J_{n_k+1} J_{n_k}$, then $(E_{n_j})$ satisfies the same 
hypotheses as $(E_\ell)$, replacing ``$\lambda_k$'' by $\tilde\lambda_k 
= \prod_{i=n_k}^{n_{k+1}-1} \lambda_i$ for all $k$.) 

Now by a ``small perturbation'' argument, we may also assume that there 
are $\ell_1<\ell_2<\cdots$ so that $E_k\subset \M_{\ell_k} \defeq Z_k$ 
for all $k$. 
Now let $Z= (Z_1\oplus Z_2\oplus\cdots)_{c_0}$ and $F\subset (Z_1\oplus 
Z_2\oplus\cdots)_{\ell^\infty}$ be defined by 
$$F= \Big\{\{e,J_1e,J_2J_1e,\ldots\} : e\in E_1\Big\}\ .
\tag 163$$
Then setting $Y= Z+F\subset (Z_1\oplus Z_2\oplus\cdots)_{\ell^\infty}$, 
$Y$ is a 1-exact operator space (since $\lambda_k\to1$). 

Now results of E.~Kirchberg and standard techniques yield  that $Y$ is 
1-locally reflexive. 
Indeed, a standard argument yields that any 1-exact operator space embeds 
in a nuclear operator space; the results in \cite{Ki2} yield in turn 
that nuclear operator spaces are 1-locally reflexive. 
(See also the last paragraph of \cite{KR}.) 
Thus by 
Lemma~3.9 of \cite{Ro2} (and the remark following it), $Z$ is 
4-completely complemented in $Y$. 
On the other hand, the argument for Lemma~4.5 yields that if $P:Y\to Z$ 
is a linear projection, then $\|P\|_{\cb}\ge C$, a contradiction.\qed
\enddemo 

We next show that $\bK_0$ (and hence $\bK$) {\it fails\/} to admit 
completely bounded extensions from certain subspaces of particular 
separable locally reflexive operator spaces.

\proclaim{Proposition 4.12} 
There exists an operator space $\tilde Y$ which is separable $1$-locally 
reflexive, a closed linear subspace $\tilde X$, and a completely bounded map 
$T:\tilde X\to \bK_0$ so that $T$ has no completely bounded extension to 
$\tilde Y$.
\endproclaim

This result follows from our work above, known results, and the following 
elementary tool. 

\proclaim{Lemma 4.13} 
Let $X,Y$ and $\tilde Y$ be operator spaces with $X\subset Y$, and let 
$q:\tilde Y\to Y$ be a complete metric surjection; set $\tilde X = 
q^{-1}(X)$ and let $T= q|\tilde X$. 
Then if $T$ has a completely bounded (resp. bounded) extension 
$\tilde T:\tilde Y\to X$, $X$ is completely complemented (resp. complemented) 
in $Y$. 
\endproclaim

\demo{Proof} 
Let $W= \ker q$; then 
$$W\subset X\ .
\tag 164$$
Now suppose $\tilde T$ is a completely bounded (resp. bounded) extension, and 
let $\Pi :\tilde Y\to \tilde Y/\tilde W$ be the quotient map and 
$S:\tilde Y/W\to Y$ the canonical complete surjective isomorphism so that 
$$q= S\Pi\ .
\tag 165$$

By (164), we may define a map $U:\tilde Y/W\to X$ by 
$$\tilde T= U\Pi\ .
\tag 166$$

Indeed, for $f\in \tilde Y$, set $U(\Pi f) = \tilde T(f)$. 
If $f\in W$, $f\in X$, hence $\tilde T(f) = T(f) = q(f) =0$; this shows 
$U$ is well defined, and we also obtain that $U$ is completely bounded 
(resp. bounded) with $\|U\|_{\cb} = \|\tilde T\|_{\cb}$ 
(resp. $\|U\| = \|\tilde T\|$). 

Now define $P:Y\to X$ by 
$$P = US^{-1}\ .
\tag 167$$

Since $T$ is a surjective quotient map from $\tilde X$ into $X$, if we 
let $x\in X$ and choose $\tilde x\in\tilde X$ with $T\tilde x =x$, we have 
that 
$$\align
P(x) &= US^{-1}q(\tilde x) = U\Pi (\tilde x)\ \text{ by (165)} \tag 168\\
& = \tilde T (\tilde x)\ \text{ by  (166)}\\
& = T(\tilde x) =X\ .
\endalign$$
Thus $P$ is a completely bounded (resp. bounded) surjective projection.\qed
\enddemo

Proposition 4.12 follows immediately from Corollary 4.3 and the next result. 

\proclaim{Proposition 4.14} 
There exists a $1$-locally reflexive separable operator space $\tilde Y$ 
with the following property: 
Given a separable operator space $X$, if 
every completely bounded (resp. bounded) linear map from a subspace 
$\tilde X$ of $\tilde Y$ to $X$ admits a completely bounded (resp. bounded) 
linear extension to $\tilde Y$, then $X$ has the {\rm CSEP} (resp the 
{\rm MSEP}).
\endproclaim

\demo{Proof} 
Let $\tilde Y = C_1$ or $(\oplus C^n_1)_{\ell^1}$, where $C_1$ is the space 
of trace-class operators (resp. $C^n_1$ is the $n$-dimensional trace-class), 
endowed with its dual structure via $C_1 = \bK^*$ 
(resp. $(\oplus C^n_1 )_{\ell^1} = \bK_0^*)$. 
A remarkable result of M.~Junge yields that $\tilde Y$ is
1-locally reflexive (\cite{J}; see also \cite{EJR} and \cite{JM}). 
But {\it every\/} separable operator space is completely isometric to 
a quotient space of $\tilde Y$ \cite{B}. 
Proposition~4.14 now  follows immediately from Lemma~4.13.\qed
\enddemo

\remark{Remark}
We do not know if $\tilde Y$ in Proposition 4.12 may be chosen so that 
$\tilde Y^*$ is separable, or so that $\tilde Y^*$ has the CMAP. 
\endremark

\enddocument